\documentclass[12pt]{article}

\usepackage{graphicx}
\usepackage{amsmath,amssymb,amsthm}
\usepackage{mathtools}
\usepackage{here} 

\usepackage[all]{xy}
\newenvironment{spmatrix}{\left ( \begin{smallmatrix}} {\end{smallmatrix}\right )}
\newenvironment{sbmatrix}{\left [ \begin{smallmatrix}} {\end{smallmatrix}\right ]}

\newtheorem{theorem}{Theorem}[section]
\newtheorem{proposition}[theorem]{Proposition}
\newtheorem{lemma}[theorem]{Lemma}
\theoremstyle{definition}

\newtheorem{alg}[theorem]{Algorithm}{\bf}{\rm}
\newtheorem{definition}[theorem]{Definition}

\newtheorem{remark}[theorem]{Remark}
\newtheorem*{acknowledgments}{Acknowledgments}
\numberwithin{equation}{section}

\newtheorem*{DAS}{Data Availability Statement}{\small \bf}{\rm \small}

\let\olditemize\itemize
\renewcommand{\itemize}{
\olditemize
\setlength{\itemsep}{1.2pt}
\setlength{\parskip}{3pt}
\setlength{\parsep}{0pt}
}
\let\oldenumerate\enumerate
\renewcommand{\enumerate}{
\oldenumerate
\setlength{\itemsep}{1.2pt}
\setlength{\parskip}{3pt}
\setlength{\parsep}{0pt}
}
\renewcommand{\labelenumi}{(\arabic{enumi})}

\begin{document}
\title{Calculation of Veech groups and Galois invariants of general origamis}
  
 \author{Shun Kumagai
}{}



\date{}
  
\maketitle
\renewcommand{\thefootnote}{\fnsymbol{footnote}}
\footnote[0]{\it Running head: \rm Veech groups and Galois invariants of general origamis}
\footnote[0]{\rm This work was supported by JSPS KAKENHI Grant Number JP21J12260. }



\begin{abstract}

Nontrivial examples of Teichm\"uller curves have been studied systematically with notions of combinatorics invariant under affine homeomorphisms.  
An origami (square-tiled surface) induces a Teichm\"uller curve for which the absolute Galois group acts on the embedded curve in the moduli space. 
In this paper, we study general origamis not admitting pure half-translation structure. 
Such a flat surface is given by a cut-and-paste construction from origami that is a translation surface. 
We present an algorithm for the simultaneous calculation of the Veech groups of origamis of given degree. 
We have calculated the equivalence classes, the $PSL(2,\mathbb{Z})$-orbits, and some Galois invariants for all the patterns of origamis of degree $d\leq 7$. 
\end{abstract}

\section{Introduction}
\label{intro}
In this paper, we consider flat surfaces obtained by gluing edges of finite copies of Euclidian unit squares. 
This sort of flat surface is called an \textit{origami} \cite{L}. 
The main result of this paper is a series of algorithms for the simultaneous calculation of the Veech groups of origamis of given degree and calculation results on origamis of degree $d\leq7$. 
In the algorithms, we enumerate all the origamis of given degree according to a cut-and-paste construction, classify them into equivalence classes, and obtain the Veech groups by calculating the $PSL(2,\mathbb{Z})$-orbit decomposition.  
The calculation was carried out using the computer resource offered under the category of General Projects by Cyber Science Center, Tohoku University.  
We note that the calculation for $d>7$ takes more than one month. 

Different from the previous implementations on origamis \cite{KOL, DH}, our result contains the raw data \cite{GitK} of all the combinatorial representation in the equivalence classes of origamis. 
The cut-and-paste construction gives many flat surfaces from a single translation surface, and we can deal with the more variations of surfaces for each degree. 
It is observed in terms of Galois invariants such as strata in the calculation result. 
The construction can be also applied to the cases of covering translation surfaces and translation surfaces with lattice Veech groups.   
Furthermore, since the cut-and-paste process preserves the holonomy vectors of all the saddle connections,  we can expect to associate with the Edwards-Sanderson-Schmidt method \cite{ESS}. 

The theory of the Teichm\"uller spaces (since the 1940s, Teichm\"uller) brought great contributions to the study of the moduli of Riemann surfaces. 
Holomorphic quadratic differentials on Riemann surfaces play a significant role in the Teichm\"uller theory. 
One important aspect is that they appear to represent the extremal deformation in the Teichm\"uller class of a quasiconformal deformation in the sense of Teichm\"uller's existence and uniqueness theorem. 
Extremal deformations are described as the locally-affine homeomorphisms between \textit{flat surfaces}, which are surfaces equipped with $(\mathbb{Z}_2\times \mathbb{R}^2)$-manifold structures, each of which is induced from the natural coordinates of a holomorphic quadratic differential. 

The orbit of a flat surface under the affine deformation action of the group $PSL(2,\mathbb{R})$ is embedded as a disk in the Teichm\"uller space with respect to the complex structure and the Teichm\"uller metric (and thus the Kobayashi metric). 
The \textit{Veech group} of a flat surface represents the action of the mapping classes on the $PSL(2,\mathbb{R})$-orbit. 
It is a discrete $PSL(2,\mathbb{R})$-subgroup \cite{V}, and its Fuchsian model corresponds to the projected image of the $PSL(2,\mathbb{R})$-orbit in the moduli space. 
Some claims on flat surfaces can be reduced to \textit{translation surfaces} through the canonical double coverings. 
Consequently the Veech group of a flat surface is a subgroup of the Veech group of its canonical double, but the converse is nontrivial.  

Examples of Veech groups have been studied systematically with notions of invariants of affine homeomorphisms. 
For a covering flat surface such as abelian origami (i.e.\ an origami that is a translation surface), one obtains the Veech group as the stabilizer of a graph structure under a combinatorial group action \cite{S1, Sh, K2}. 
Invariants in terms of the geometry of the natural flat metric on a flat surface enable us to describe the membership criterion for the Veech group of a general flat surface and evaluate the $PSL(2,\mathbb{R})$-orbit \cite{B, Mu, ESS}. 
Unlike most previous researches on Veech groups, this paper focuses on flat surfaces, not translation surfaces. 

Origamis are studied in the context of the Galois action on combinatorial objects as well as \textit{dessins d'enfants} \cite{JW}. 
The $PSL(2,\mathbb{R})$-orbit of an origami projects to an arithmetic algebraic curve embedded in the moduli space, for which the two Galois actions are compatible: for the curve itself and the moduli represented by each point of the curve. 
A crucial result is given by M\"oller \cite{M}, and an overview of this ground is described in \cite{L, HS1}. 
Examples of nontrivial Galois conjugacy of origamis are given by using `M-origamis' \cite{M,Nis} constructed from dessins where any nontrivial Galois conjugacy is inherited. 
Considering Veech groups of origamis is also meaningful in the following sense: 
It is known that a flat surface admits two directional cylinder decompositions if the Veech group is a lattice \cite{B}. 
Then, the surface is up to affine deformation uniquely represented by an origami with its square cells replaced by rectangles of possibly distinct moduli \cite{K2}. 
The membership criterion for the Veech group is stated using origamis as Proposition \ref{VG_determined}.


The structure of this paper is as follows: 
In section \ref{sec:2}, we introduce basic discussions on flat surfaces and origamis. 
The concept of the main algorithms is based on the observation stated in Remark \ref{importantrem}. 
In section \ref{sec:3}, we discuss the cut-and-paste construction of an origami from an abelian origami. 
In section \ref{sec:4}, we present the concrete steps of the main algorithms. 
In section \ref{sec:5}, we present the calculation result and consider the Galois invariants of origamis of degree $d\leq 7$. 

\section{Preliminaries}
\label{sec:2}
In this section, we introduce basic notations and results on flat surfaces for the main discussions. 
\subsection{flat surface}
\label{sec:2.1}

Let $R$ be a Riemann surface of type $(g,n)$ with $2g-2+n>0$.

\begin{definition}
A \textit{holomorphic quadratic differential} $\phi$ on $R$ is a tensor on $R$ whose restriction to each chart $(U,z)$ on $R$ is of the form $\phi_U(z)dz^2$ where $\phi_U:U\rightarrow \hat{\mathbb{C}}$ is holomorphic. 
A pair $(R,\phi)$ is called a \textit{flat surface}. 
We say that the set $\mathrm{Sing}(R,\phi)$ of marked points of $R$ and zeros and poles of $\phi$ is the set of the \textit{singularities} of $(R,\phi)$. 
\end{definition}

Let $p_0\in R^*= R\setminus\mathrm{Sing}(R,\phi)$ and $(U,z)$ be a chart around $p_0$. 
Then, $\phi$ defines a natural coordinate ($\phi$\textit{-coordinate}) on $U$ by 
\begin{align}
\zeta_\phi(p)=\int^p_{p_0} \sqrt{\phi_U(z)}dz\text{,\ \ $p\in U$},\end{align}
for which $\phi=(d\zeta_\phi)^2$ holds. 
The $\phi$-coordinates give an atlas ${A}_\phi$ on $R^*$ any of whose transition map is a \textit{half-translation} $\zeta \mapsto \pm \zeta +c\ (c\in \mathbb{C})$. 
Such a structure, an atlas any of whose coordinate transformation is a half-translation is called a \textit{flat structure}. 
The atlas ${A}_\phi$ extends to each singularity $p_0\in\mathrm{Sing}(R,\phi)$ of order $m$, with local representation
\begin{equation}\label{cone}
	\zeta_\phi(p)=\int^p_{p_0} \sqrt{z^m}dz=z(p)^{\frac{m}{2}+1},\text{\ \ $p\in U\setminus\{ p_0\}$}, 
\end{equation}
for a suitable chart $(U,z)$ around $p_0$.

\begin{definition}Let $(R,\phi)$ and $(S,\psi)$ be flat surfaces of genus $g$. 
\begin{enumerate}
\item For $A=\begin{spmatrix}a&b\\ c&d\end{spmatrix}\in GL(2,\mathbb{R})$, we denote by $[A]=\begin{sbmatrix}a&b\\ c&d\end{sbmatrix}$ the quotient class of $A$ in $PSL(2,\mathbb{R})$. 
We define $f_A:\mathbb{C}\rightarrow\mathbb{C}$ by
		\begin{equation}
		f_A(\xi+{i}\eta)= (a\xi+c\eta)+{i}(b\xi+d\eta), \ \ \xi,\eta\in\mathbb{R}. 
		\end{equation}
	\item We say that a branched covering $f:(R,\phi)\rightarrow (S,\psi)$ is \textit{locally-affine} if there exist finite 
{  sets  $\Sigma_R$, $\Sigma_S$}
 such that 
{  $\mathrm{Sing}(R,\phi)\subset \Sigma_R\subset R$, $\mathrm{Sing}(S,\psi)\subset \Sigma_S\subset S$,  and}
$f$ is restricted to a covering $f:R\setminus\Sigma_R\rightarrow S\setminus\Sigma_S$ that is locally represented by $z \mapsto f_A(z) + c$  for some $A\in SL(2,\mathbb{R})$ and $c\in \mathbb{C}$ with respect to natural coordinates of $\phi$ and $\psi$. 
A locally-affine biholomorphism is called a \textit{half-translation map}, and we say that two flat surfaces are \textit{equivalent} if there exists a half-translation map between them.
	\item For a locally-affine covering $f:(R,\phi)\rightarrow (S,\psi)$, the local derivative $A$ on $R\setminus\Sigma_R$ is constant up to a factor 
	$\{\pm 1\}$. 
	We call $D(f):=[A]\in PSL(2,\mathbb{R})$ the \textit{derivative} of $f$. 
	\item The group of locally-affine, orientation-preserving self-homeomorphisms of $(R,\phi)$ is called the \textit{affine group} of $(R,\phi)$ and denoted by $\mathrm{Aff}^+(R,\phi)$. 
	The group $\Gamma (R,\phi)
$ of the derivatives of elements in $ \mathrm{Aff}^+(R,\phi)$  is called the 
	\textit{Veech group} of $(R,\phi)$. 
\end{enumerate}
\end{definition}

For a flat surface $(R,\phi)$, the family of affine homeomorphisms on $(R,\phi)$ forms a disk $\Delta(R,\phi)$ isometrically embedded in the Teichm\"uller space of $R$ by Teichm\"uller's theorem. 
The stabilizer of $\Delta(R,\phi)$ in the Teichm\"uller-modular group is known to be the affine group $\mathrm{Aff}^+(R,\phi)$.
Furthermore, the action $\mathrm{Aff}^+(R,\phi)$ on $\Delta(R,\phi)$ is given by the Veech group $\Gamma(R,\phi)$ acting on the unit disk by M\"obius transformations.
(See \cite{GL} for instance.)

It is first observed by Veech \cite{V} that a Veech group is a discrete group. 
In particular, the projected image $C(R,\phi)$ of $\Delta(R,\phi)$ in the  moduli space is an orbifold isomorphic to 
{  the mirror image of}
$\mathbb{H}/\Gamma (R,\phi)$.
If $\Gamma (R,\phi)$ has finite covolume, $C(R,\phi)$ is an algebraic curve called the \textit{Teichm\"uller curve}. 

\begin{definition}
We say that a flat surface $(R,\phi)$ is \textit{abelian} if $\phi$ becomes the square of an abelian differential on $R$ and otherwise \textit{non-abelian}. 
\end{definition}
In this paper, we consider integrable flat surfaces, equivalently flat surfaces having singularities of order no less than $-1$. 
The space $\mathcal{Q}_g$ of equivalence classes of integrable flat surfaces of genus $g$ 
is naturally strartified as 
\begin{equation} \label{stratification}
	\mathcal{Q}_g=\left({\bigsqcup} \mathcal{A}_g(2m_1,\ldots 2m_k) \right)\sqcup\left({\bigsqcup} \mathcal{Q}_g(m_1',\ldots m_{k'}') \right),\end{equation}
where $\mathcal{A}_g(2m_1,\ldots 2m_k)$ (resp.\ $\mathcal{Q}_g(m_1',\ldots m_{k'}')$) is the {stratum} of abelian (resp.\ non-abelian) flat surfaces of genus $g$ with precisely $k$ (resp.\ $k'$) singularities of orders $2m_1\leq\cdots\leq 2m_k$ (resp.\ $m_1'\leq\cdots\leq m_{k'}'$). 
It follows from Riemann-Roch theorem that the orders of singularities of a flat surface of genus $g$ sum up to $4g-4$. 
Thus the indices in the stratification (\ref{stratification}) run through integers such that 
\begin{align}\label{strata}
	m_1,\ldots ,m_k\in \mathbb{N}, &&m_1',\ldots ,m_{k'}'\in\{-1\}\cup\mathbb{N}, &&\sum_{j=1}^{k}2m_j=\sum_{j=1}^{k'}m_j'=4g-4.
\end{align}
Kontzevich-Zorich \cite{KZ} and Lanneau \cite{Lanneau} showed that strata $\mathcal{A}_g(2m_1,\ldots 2m_k)$, $\mathcal{Q}_g(m_1',\ldots m_{k'}')$ with (\ref{strata}) are noempty except for few exceptional cases. 
\begin{definition}\label{canonical_double}
The \textit{canonical double covering} of a flat surface $(R,\phi)$ is the double covering $\pi_{\phi}:\hat{R}\rightarrow R$ obtained from the continuation of branches of locally defined abelian differential $\sqrt{\phi}$ on $R$. 
\end{definition}
By definition, the canonical double covering $\pi_{\phi}:\hat{R}\rightarrow R$ induces a holomorphic abelian differential $\hat{\phi}=\pi_\phi^*\sqrt{\phi}$ on $\hat{R}$. 
As a consequence of \cite[Theorem 3.6]{GJ}, a canonical double covering satisfies the following property. 
\begin{proposition}\label{lifting}
Let $(R,\phi)$ be a flat surface, $(S,\psi)$ be an abelian flat surface, and $f:(S,\psi)\rightarrow (R,\phi)$ be a locally-affine covering. 
Then $f$ is the composition of a locally-affine covering $\hat{f}:(S,\psi)\rightarrow (\hat{R},\hat{\phi})$ and the canonical double covering $\pi_{\phi}:\hat{R}\rightarrow R$. 
\end{proposition}
\begin{remark}\label{poly_double}
By continuing natural coordinates, every flat surface is represented by finitely many Euclidian polygons $(C_i)_{i\in I}$ with their edges glued by half-translations. 
On this representation, the canonical double covering is obtained as follows (see Figure \ref{figure3}):
 \begin{enumerate}
	\item For each $i\in I$, take a copy $C^+_i$ and a half-rotated copy $C^-_i$ of the polygon $C_i$. 
	Denote by $e_+$ (resp.\ $e_-$) the edge of the polygon $C^+_i$ (resp.\ $C^-_i$) corresponding to an edge $e$ of the polygon $C_i$. 
	\item If an edge $e$ of a polygon $C_i$ is glued with an edge $e'$ of a polygon $C_{i'}$ by transition map of the form $z\mapsto z+c$ (resp.\ $-z+c$), glue $e_+$ with $e'_+$ and $e_-$ with $e'_-$ (resp.\ $e_+$ with $e'_-$ and $e_-$ with $e'_+$) by translation. 
  \end{enumerate}
\begin{figure}[htbp]
\begin{center}
  \includegraphics[width=130mm]{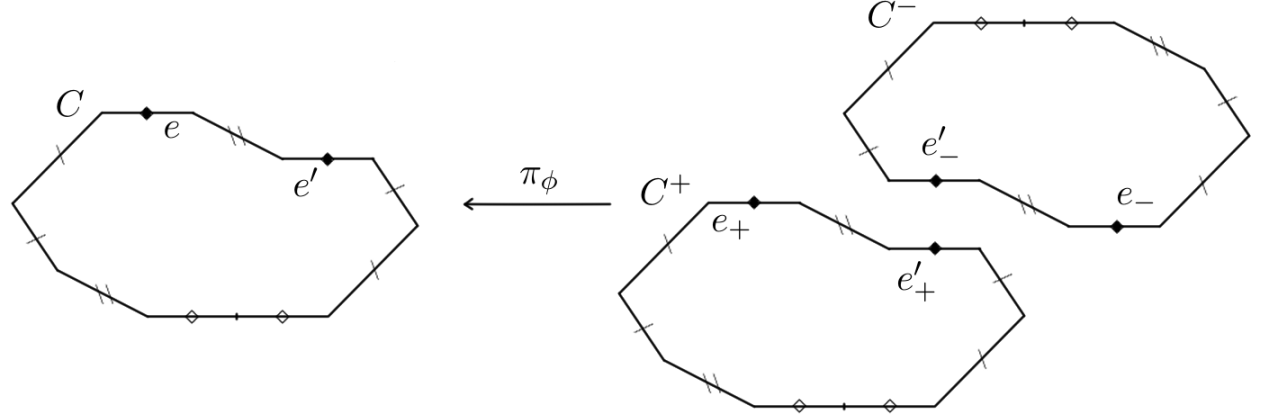}
  \caption{The canonical double covering of a flat surface represented by a Euclidian polygon: edges of $C$ with the same marks are glued. }
\label{figure3}
\end{center}
\end{figure}
\end{remark}

\subsection{$\phi$\textit{-metric}}
\label{sec:2.2}

Let $(R,\phi)$ be a flat surface of genus $g$. 
The Euclidian metric lifts via $\phi$-coordinates to a flat metric 
on $R$, called the $\phi$\textit{-metric}. 
A geodesic of $\phi$-metric is called a $\phi$\textit{-geodesic}. 
Via the $\phi$-coordinates, a 
$\phi${-geodesic} 
{  without singularities} 
is locally a line segment on the plane whose \textit{direction} is uniquely determined in $\mathbb{R}/\pi\mathbb{Z}$. 

\begin{definition}Let $(R,\phi)$ be a flat surface of genus $g$. 
\begin{enumerate}
	\item The $\phi$\textit{-cylinder} generated by a $\phi$-geodesic $\gamma$ is the union of all $\phi$-geodesics parallel (with same direction) and free homotopic to $\gamma$. 
We define the direction of a 
{  $\phi$-cylinder} 
by the one of 
{  a generator of it.} 
	\item $\theta\in \mathbb{R}/\pi\mathbb{Z}$ is called \textit{Jenkins-Strebel direction} of $(R,\phi)$ if almost every point in $R$ lies on some closed $\phi$-geodesic in the direction $\theta$. 
	Let $JS(R,\phi)$ denote the set of Jenkins-Strebel directions of $(R,\phi)$. 
\end{enumerate}
\end{definition}
{  Let $f$ be a locally-affine map on a flat surface $(R,\phi)$. Then the map $f$ maps any line segment in the direction $\theta\in \mathbb{R}/\pi\mathbb{Z}$ to a line segment in the direction $A\theta:=\mathrm{arg}(f_{D(f)}(e^{\sqrt{-1}\theta}))$.} 
If $(R,\phi)$ admits two Jenkins-Strebel directions $\theta_1, \theta_2\in JS(R,\phi)$, then the surface is decomposed into parallelograms which are intersections of cylinders in the directions $\theta_1, \theta_2$. 
By continuing local segments in each directions, we may define a combinatorial structure of such a decomposition similar to an origami (see section \ref{sec:origami}). 
This structure is uniquely determined by the class of $(R,\phi)$, and we may determine the existence of an affine map of prescribed derivative by comparing decompositions as follows. 
\begin{proposition}[\cite{K2}]\label{VG_determined}
Let $(R,\phi)$ be a flat surface with two distinct finite Jenkins-Strebel directions $\theta_1,\theta_2\in JS(R,\phi)$. 
A matrix $A\in PSL(2,\mathbb{R})$ belongs to $\Gamma (R,\phi)$ if and only if the following holds. 
\begin{enumerate}
\item $A\theta_1,A\theta_2$ 
{ belong}
 to $ JS(R,\phi)$. 
\item The origamis with compatible moduli lists given by the decomposition of $(R,\phi)$ in $(\theta_1,\theta_2)$ and $A(\theta_1,\theta_2):=(A\theta_1,A\theta_2)$ are equivalent. 
That is, there exists a permutation that conjugates two origamis with moduli lists. 
\end{enumerate}
\end{proposition}

To calculate the Veech group of a flat surface $(R,\phi)$, it is good to consider the action of $\mathrm{Stab}_{PSL(2,\mathbb{R})}JS(R,\phi)$ on the set of parallelogram decompositions of $(R,\phi)$ defined by taking the decomposition in $A^{-1}(\theta_1, \theta_2)$ for each $\theta_1, \theta_2\in JS(R,\phi)$ and $A\in PSL(2,\mathbb{R})$. 
The image under $A$ determines the class of $(R,f_A^*\phi)$, and the stabilizer of a decomposition is the Veech group of $(R,\phi)$.

\subsection{origami}
\label{sec:origami}

\begin{definition}
An \textit{origami} of degree $d$ is a flat surface obtained by gluing $d$ Euclidian unit squares at edges equipped with the flat structure induced from the natural coordinates of squares.
\end{definition}
\begin{figure}[htbp]
\begin{center}
  \includegraphics[width=100mm]{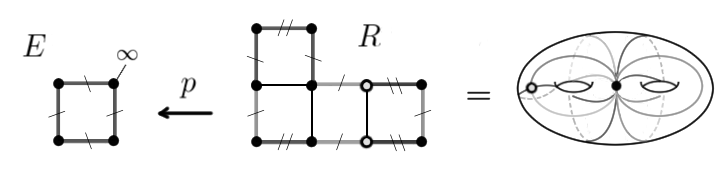}
  \caption{An origami of degree $4$: edges with the same marks are glued. }
\label{Figure1}
\end{center}
\end{figure}
In this paper, we shall mark all the corner points of an origami by treating nonsingular corner points as singularities of order $0$. 
Figure\ \ref{Figure1} shows an example of origami of degree $4$ that belongs to the stratum $\mathcal{A}_2(0,4)$. 
For an abelian origami, a cell-to-cell correspondence defines a locally-affine covering $p:R\rightarrow E$ of the unit square torus $(E=\mathbb{C}/(\mathbb{Z}+i\mathbb{Z}),dz^2)$ branched at most over one point $\infty\in E$.  
It induces the following characterizations of an abelian origami. (See \cite{HS1} for details.)


\begin{lemma}\label{origami}
An { abelian} origami of degree $d$ is up to equivalence (in brackets) uniquely determined by each of the following. 
\begin{enumerate}\renewcommand{\labelenumi}{\rm(\alph{enumi})}
	\item {A connected, oriented graph $(\mathcal{V, E})$ with $|\mathcal{V}|=d$ such that each vertex has 
{  precisely two incoming edges and two outgoing edges, one of each of which is labeled with $x$ and the other with $y$.}
[up to  equivalence of the natural filling graph embedding]. }
	\item{A $d$-fold branched covering $p:R\rightarrow E$ of a torus $E$ branched at most over one point $\infty\in E$ [up to covering equivalence]. }
	\item{A pair $(x,y)$ of two elements in the symmetric group $\mathfrak{S}_d$ generating a transitive permutation group $G$ [up to conjugation in $\mathfrak{S}_d$]. }
	\item{A subgroup $H$ of the free group $F_2$ of index $d$ [up to conjugation in $F_2$]. }
\end{enumerate}
\end{lemma}

The two generators $x,y\in G$ correspond to the monodromies along the horizontal and vertical cylinder curves on $E$, respectively. 
The permutation $z=xyx^{-1}y^{-1}\in \mathfrak{S}_d$ is the monodromy around the corner point of square cells.  
In particular, an abelian origami belongs to the stratum $\mathcal{A}_g(2m_1,\ldots,2m_n)$ if and only if $z$ consists precisely of $n$ cycles of lengths $m_1+1,\ldots,m_n+1$.

\begin{remark}\label{importantrem}\ 
\begin{enumerate}
\item
For an origami $(R,\phi)$, the set of Jenkins-Strebel directions is $JS(R,\phi)=PSL(2,\mathbb{Z})\cdot (0+\pi\mathbb{Z})\subset \mathbb{R}/\pi\mathbb{Z}$. 
By Proposition \ref{VG_determined}, the Veech group of an origami is a subgroup of $PSL(2,\mathbb{Z})$.
\item
Let $\bar{\Omega}_d$ be the set of all the classes of origamis of degree $d\in\mathbb{N}$.
$PSL(2,\mathbb{Z})$ acts on $\bar{\Omega}_d$ by $A\in PSL(2,\mathbb{Z}):\mathcal{O}\mapsto \mathcal{O}_A$ where $\mathcal{O}_A$ is the origami obtained as the parallelogram decomposition of $\mathcal{O}$ in $A^{-1}(0,\frac{\pi}{2})$. 
By Proposition \ref{VG_determined}, the Veech group of an origami $\mathcal{O}$ is the stabilizer of the equivalence class of $\mathcal{O}$ under this action. 
We may apply the Reidemeister-Schreier method \cite{MKS} to the result of Algorithm \ref{alg_comp} to obtain the list of generators and the list of representatives of the Veech group of each origami. 
\item Schmith\"usen \cite{S1} showed that the action of $PSL(2,\mathbb{Z})$ on the set of abelian origamis is desctribed by the automorphisms $\gamma_T,\gamma_S$ on $F_2=\langle x,y\rangle$ defined by
\begin{align}
  \gamma_T(x,y)=(x,xy),\ \ \ \ \gamma_S(x,y)=(y,x^{{-1}}), 
\end{align}
where $[T]=\begin{sbmatrix}1&1\\ 0&1\end{sbmatrix}$ and $[S]=\begin{sbmatrix}0&-1\\ 1&0\end{sbmatrix}$ are matrices generating $PSL(2,\mathbb{Z})$. 
\end{enumerate}
\end{remark}

Let $\mathcal{O}$ be an arbitrary origami. 
As the Veech group $\Gamma(\mathcal{O})$ is a finite index $PSL(2,\mathbb{Z})$-subgroup, it induces a Teichm\"uller curve $C(\mathcal{O})$ as a Bely\u{\i} covering $\beta_{C(\mathcal{O})}:\mathbb{H}/\Gamma(\mathcal{O})\rightarrow \mathbb{H}/PSL(2,\mathbb{Z})$ branched at most over the three points. 
The origami $\mathcal{O}$ also admits a Bely\u{\i} covering $\beta_\mathcal{O}$ induced from the rational map $\beta(x,y)=4x^2$ on 
{ the elliptic curve $y=4x^3-x$ that represents the unit square torus.}
Both are algebraic curves defined over $\bar{\mathbb{Q}}$ by Bely\u{\i} theorem. 
M\"oller \cite{M} proved the compatibility of the action of the absolute Galois group $G_\mathbb{Q}=\mathrm{Gal}(\bar{\mathbb{Q}}/\mathbb{Q})$ for the curve $C(\mathcal{O})$ embedded in the moduli space. 
That is, for any $\sigma\in G_{\mathbb{Q}}$, the curve $C(\mathcal{O})^\sigma$ coincides with the Teichm\"uller curve $C(\mathcal{O}^\sigma)$ induced from the origami $\mathcal{O}^\sigma$ and every point $\mathcal{O}_t\in C(\mathcal{O})$ that corresponds to a curve defined over $\bar{\mathbb{Q}}$ is conjugated to a point $\mathcal{O}_t^\sigma\in C(\mathcal{O})^\sigma$. 

The \textit{valency list} of a Bely\u{\i} covering $\beta:C\rightarrow \hat{\mathbb{C}}$ is the list
$(k_1^{0},\ldots ,k_{l_0}^{0}\mid k_1^{1},\ldots ,k_{l_1}^{1} \mid k_1^{\infty},\ldots k_{\infty}^{\infty})$ of all the ramification indices over the three branched points. 
The genus of $C$, the degree and the valency list of $\beta$ are classically known as $G_\mathbb{Q}$-invariants. 
It is classically known that the action of $G_\mathbb{Q}$ preserves local behavior of a branched covering. 
With M\"oller's work, 
we can see that the valency list of $\beta_{C(\mathcal{O})}$, the genus of $C(\mathcal{O})$, and the stratum of $\mathcal{O}$ are  $G_\mathbb{Q}$-invariants of an origami $\mathcal{O}$.
\medskip\medskip
\section{Cut-and-paste construction}
\label{sec:3}
In this section, we discuss a cut-and-paste construction of origamis and prove formulae for specifying the equivalence classes and the strata. 

\begin{definition}\label{formula_xye}  Let $I_d:=\{1,\ldots ,d\}$ be the set of $d$ indices and $\bar{I}_d{:=\{\pm1\}\times I_d}=\{\pm1,\ldots ,\pm d\}$ be its double. 
  Let $\mathcal{E}_{d}:=\{\varepsilon{ \in}\{\pm1\}^{\bar{I}_d}\mid\varepsilon(- i )={ -} \varepsilon( i ),\  i \in\bar{I}_d \}$ be the set of {  odd}  functions on $\bar{I}_d$. 
Let $\bar{\mathfrak{S}}_d:=\mathrm{Sym}(\bar{I}_d)$ naturally embed the symmetric group $\mathfrak{S}_d$  { in the set ${ \bar{\mathfrak{S}}^{\mathrm{odd}}_{d}}$ of odd functions in $\bar{\mathfrak{S}}_d$}. 
  For each $x\in \mathfrak{S}_d$ and $\varepsilon\in\{\pm1\}^{\bar{I}_d}$, let $x^\varepsilon$ denote the mapping on $\bar{I}_d$ defined by 
\begin{equation}
x^\varepsilon( i ) {:=  (x^{\varepsilon(i)})( i ) }
\text{ for each } i  \in \bar{I}_d.
\end{equation}
\end{definition}

\begin{figure}[htbp]
  \begin{center}
  \includegraphics[width=120mm]{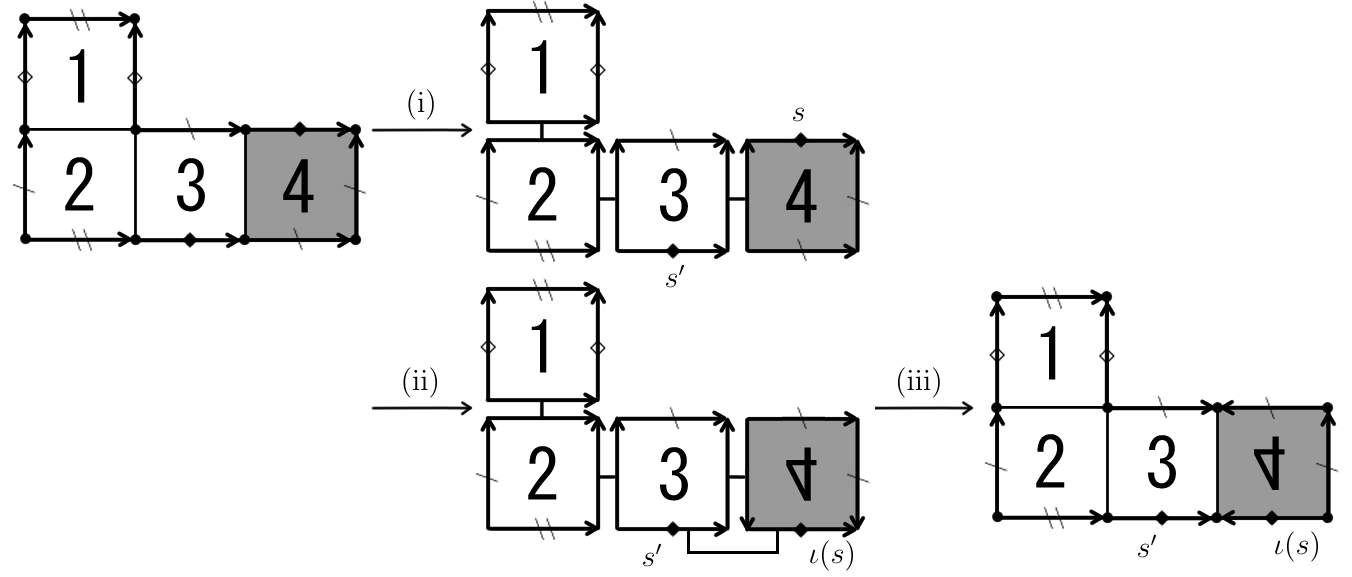}
    \caption{The origami given by $(x,y, \varepsilon)$ where $x=(1)(2\ 3\ 4)$, $y=(1\ 2)(3\ 4)$, $\varepsilon=(+,+,+,-)$: edges with the same marks are glued. 
{ 
In the step (i), the upper side $s$ of the 4th cell is paired with the lower side $s'$ of the 3rd cell. 
In the step (ii), the vertical reflection $\iota$ is applied to the 4th cell. 
In the step (iii), the sides $\iota(s)$ and $s'$ are glued by a half-translation. }}
  \label{construction}    
  \end{center}
  \end{figure}

\begin{definition}\label{xye_def}Let
$\Omega _d:={\mathfrak{S}_d}\times {\mathfrak{S}_d}$ be the set of (possibly disconnected) abelian origamis of degree $d$, $\Omega^0_{2d}:=\{\mathcal{O}\in\Omega_{2d}\mid$ there exists an origami whose canonical double covering surface is $\mathcal{O}\}$, and $\bar{\Omega}_d:=\Omega _d\times \mathcal{E}_d$.
For each $\mathcal{O}=(x,y,\varepsilon)\in\bar{\Omega}_d$ and $i  \in \bar{I}_d$, define  
\begin{equation}\label{xye_eq}
  \left\{
  \begin{array}{rl}
     \mathbf{x}_\mathcal{O}( i )&=x^{\mathrm{sign}}( i ) \\
     \mathbf{y}_\mathcal{O}( i )&=\varepsilon( i )\cdot y^{\varepsilon}( i )\cdot \varepsilon(y^{\varepsilon}( i ))\ ,
  \end{array}
  \right.
\end{equation}
{ where $\mathrm{sign}\in\mathcal{E}_d$ is the sign function.}
\end{definition}
Observe that a \textit{double symmetry} $\mathbf{w}(-\mathbf{w}(-i))=i$, $i\in \bar{I}_d$ holds for $\mathbf{w}=\mathbf{x}_\mathcal{O},\mathbf{y}_\mathcal{O}$ and thus (\ref{xye_eq}) defines two permutations $\mathbf{x}_\mathcal{O},\mathbf{y}_\mathcal{O}\in \bar{\mathfrak{S}}_d$. 
We may assign each triple $\mathcal{O}=(x,y,\varepsilon)\in\bar{\Omega}_d$ to an origami of degree $d$ whose canonical double covering is the abelian origami $(\mathbf{x}_\mathcal{O},\mathbf{y}_\mathcal{O})\in\Omega_{2d}$ in the following way (see Figure\ \ref{construction}): 
  \begin{enumerate} \setlength{\leftskip}{5pt}\renewcommand{\labelenumi}{(\arabic{enumi})}
    \item[(i)] Cut the abelian origami $(x,y)$ at all edges (with the edge-pairings remembered). 
    \item [(ii)]Apply the vertical reflection to the $i$-th cell if $\varepsilon(i)=-1$ for each $i\in I_d$. 
    \item [(iii)]Glue all paired edges in such a way that with the natural coordinates, the quadratic differential $(dz)^2$ is globally defined on the resulting surface. 
 \end{enumerate}

If we add a half-rotated copy with cells labelled by $\{-1,\ldots ,-d\}$, we can make the similar construction so that the abelian differential $dz$ is globally defined on the resulting surface. 
In the construction, every horizontal sides remain the same, and $x$ is uniquely extended to $\mathbf{x}_\mathcal{O}$ by double symmetry. 
{ Fix natural coordinates around upper and lower sides of all the cells in the step (i). 
After the step (ii), each coordinate around an upper side $s$ of the $i$-th cell is multiplied by $\varepsilon(i)$ by the vertical reflection. 
The side $s$ is glued with the lower side $s'$ of either $y^{\varepsilon(i)}(i)$-th cell or $-y^{\varepsilon(i)}(i)$-th cell, so that the corresponding transition map is a translation. 
Since each coordinate around the side $s'$ is multiplied by $\varepsilon(y^{\varepsilon(i)}(i))$, the side $s'$ belongs to the $\varepsilon(i)\cdot y^{\varepsilon(i)}(i)\cdot \varepsilon(y^{\varepsilon(i)}(i))$-th cell. 
} 
In this way, we obtain the abelian origami represented by $(\mathbf{x}_\mathcal{O},\mathbf{y}_\mathcal{O})$ as to be the canonical double covering
{ of the origami $\mathcal{O}$}
.  

Conversely, an abelian origami $(\mathbf{x,y})$ admits an involutive deck transformation $\tau$ locally represented by $z\mapsto -z$. 
We may label the cells of the origami $(\mathbf{x,y})$ in $\bar{I}_d$ so that $\tau$ acts as $i\mapsto -i$, and then $\mathbf{x,y}$ have double symmetry. 
For an abelian origami $(\mathbf{x,y})$ with double symmetry, a triple $(x,y,\varepsilon)$ is reconstructed by the following lemma.

\begin{lemma}\label{restore}
For a cyclic permutation $a={(a_1\ a_2\ \cdots\ a_m)}$, denote $a':=(-a_1 \ -\!a_2\ \cdots\ -\!a_m)^{-1}$ and $|a|:=(|a_1|\ |a_2|\cdots |a_m|)$.
  Let $\mathbf{y}\in\bar{\mathfrak{S}}_d$ be a permutation with double symmetry and fix a cycle decomposition $\mathbf{y}=c_1{c_1'}\cdots c_n{c_n'}$. 
For $j=1,\ldots ,n$, define $\varepsilon_j\in\mathcal{E}_d$ so that  
{ 
    \[
    \varepsilon_j( i )= 
    \left\{
     \begin{array}{rl}
        -1&\text{ if $c_j(- i )\neq - i $}\\
        1&\text{ otherwise}
     \end{array}
    \right. 
\text{,\ \  for each $ i \in \bar{I}_d$.}
    \]
}
  Then, $y={ \check{\mathbf{y}}}:=|c_1|\cdots |c_n|\in{ \mathfrak{S}_d}$ and $\varepsilon=\varepsilon_\mathbf{y}:=\varepsilon_1\cdots \varepsilon_n \in \mathcal{E}_d$ satisfy that 
\begin{equation}
     \mathbf{y}( i )=\varepsilon( i )\cdot y^{\varepsilon}( i )\cdot \varepsilon(y^{\varepsilon}( i ))
  \text{ for each } i  \in \bar{I}_d.
\end{equation}
\end{lemma}
\begin{proof}
  Suppose $n=1$. 
  Denote $\mathbf{y}=cc'=(a_1\ a_2\ \cdots\ a_d)(-a_1\ -\!a_2\ \ldots \ -\!a_d)^{-1}$, $y=(|a_1|\ |a_2|\cdots\ |a_d|)(-|a_1|\hspace{2pt} -\!|a_2|\cdots\hspace{2pt}-\!|a_d|)^{ -1}$, and $a_{d+1}=a_1$.
{  Since $\varepsilon(i)=-1$ if and only if $i$ belongs to the cycle $c'$, }
we have $\varepsilon (a_i)=1$ and $\mathbf{y}(a_i)=\varepsilon(|a_{i+1}|)|a_{i+1}|$, for all $i\in I_d$.
  \text{For each $i\in I_d$, }
  \begin{align}
    \mathbf{y}(y,\varepsilon)(a_i)&{ :=}\varepsilon(a_i)\cdot y^{\varepsilon(a_i)} (a_i)\cdot\varepsilon(y^{\varepsilon(a_i)} (a_i))\nonumber\\
    &=y(\mathrm{sign}(a_i)|a_i|)\cdot \varepsilon(y(\mathrm{sign}(a_i)|a_i|))\nonumber\\
    &=\mathrm{sign}(a_i)|a_{i+1}|\cdot \mathrm{sign}(a_i)\varepsilon(|a_{i+1}|)\nonumber\\
    &=\varepsilon(|a_{i+1}|)|a_{i+1}|\nonumber\\
    &=\mathbf{y}(a_i).
  \end{align}
  Thus $\mathbf{y}(y,\varepsilon){ =}\varepsilon\cdot y^\varepsilon \cdot\varepsilon(y^\varepsilon)$ equals to $\mathbf{y}$. 
  Applying this to each cycle in $\mathbf{y}=c_1{c_1'}\cdots c_n{c_n'}$, we obtain the claim for general $n$.\end{proof}

{Remark that ${ \check{\mathbf{y}}},\varepsilon_{\mathbf{y}}$ {  depend} on the way choosing cycles $c_1,\ldots ,c_n$. }
For any $\mathbf{x}\in\bar{\mathfrak{S}}_d$ with double symmetry, we may take $x={ \check{\mathbf{x}}}$ {  as in Lemma \ref{restore}} so that $\mathbf{x}(i)=x^{\mathrm{sign}}(i)$ for each $i\in \bar{I}_d$. 
In particular, 
$\mathcal{O}\mapsto (\mathbf{x}_\mathcal{O},\mathbf{y}_\mathcal{O})$ defines a 1-1 correspondence $\bar{\Omega}_d\rightarrow \Omega^0_{2d}$ up to equivalence. 

{ 
\begin{definition}
For $\sigma,\tau\in \bar{\mathfrak{S}}_d$ and $\varepsilon\in \mathcal{E}_d$ , we denote $\sigma^\#\tau:=\sigma\circ\tau\circ\sigma^{-1}\in\bar{\mathfrak{S}}_d$ and $\xi(\tau,\varepsilon):=\varepsilon\cdot(\varepsilon\circ\tau)\in\{\pm1\}^{\bar{I}_d}$. 
If $\sigma$ belongs to ${ \bar{\mathfrak{S}}^{\mathrm{odd}}_{d}}$, we define ${\delta_\sigma}{:=\mathrm{sign}(\sigma(|\cdot|))}\in\{\pm1 \}^{\bar{I}_d}$ (even function) and ${{\bar{\sigma}}}{:={\delta_\sigma}\sigma} \in \mathfrak{S}_d$. 
\end{definition}
We note that  $\mathbf{y}(y,\varepsilon)=\xi(y^{\varepsilon},{\varepsilon})\cdot y^{\varepsilon}$, $\xi(y^\varepsilon,\varepsilon)(-i)=\xi(y^{-\varepsilon},\varepsilon)(i)$, and $\xi(\tau,\lambda)\xi(\tau,\lambda')=\xi(\tau,\lambda\lambda')$ hold. 
}

  \begin{lemma}\label{isom}Let $\mathcal{O}_j=(x_j,y_j,\varepsilon_j)\in\bar{\Omega}_d\ (j=1,2)$ be two origamis. Then $\mathcal{O}_1,\mathcal{O}_2$
  are equivalent as flat surfaces if and only if there exists $\sigma={ \delta_\sigma}\bar{\sigma}\in{ \bar{\mathfrak{S}}^{\mathrm{odd}}_{d}}$
   such that the following (1)-(4) holds on $I_d$: 
    \begin{enumerate}\setlength{\leftskip}{5pt}\setlength{\itemsep}{3pt}\renewcommand{\labelenumi}{(\arabic{enumi})}
    \item ${ \delta_\sigma}={ \delta_\sigma}\circ x_1$,
    \item $x_2={{\bar{\sigma}}}^\# (x_1^{{ \delta_\sigma}})$,
    \item 
    $\xi(y_2,{ \delta_\sigma}\circ{{\bar{\sigma}}}^{-1}\cdot \varepsilon_1\circ{{\bar{\sigma}}}^{-1}\cdot\varepsilon_2)=1$, and
    \item $y_2={{\bar{\sigma}}}^\#(y_1^{{ \delta_\sigma}\cdot\varepsilon_1\cdot\varepsilon_2\circ{{\bar{\sigma}}}})$. 
    \end{enumerate}
\end{lemma}
\begin{proof}
Assume that there exists a half-translation map between $\mathcal{O}_1$ and $\mathcal{O}_2$.
By Proposition \ref{lifting}, it lifts via their canonical double coverings and induces a cell-to-cell correspondence $\sigma\in\bar{\mathfrak{S}}_{d}$ such that 
$\mathbf{x}_2(i)=\sigma^\#\mathbf{x}_1(i)$ and $\mathbf{y}_2( i)=\sigma^\#\mathbf{y}_1( i)$ for $i\in{I}_d$. 
Since $\mathbf{x}_1$ and $\mathbf{x}_2$ have double symmetry, it follows that $\mathbf{x}_2(\sigma(-i))=\mathbf{x}_2(-\sigma(i))$ for each $i\in I_d$ and thus { $\sigma\in{ \bar{\mathfrak{S}}^{\mathrm{odd}}_{d}}$}.
For $i\in{I}_d$ and $\varepsilon\in\{\pm1 \}$, we have 
the following:
\[
 \begin{array}{rlr}
    \sigma(\mathbf{x}_1)(\varepsilon i)&=({ \delta_\sigma}{{\bar{\sigma}}})(x_1^{\mathrm{sign}(\varepsilon i)}(\varepsilon i))\\[4pt]
    &=\varepsilon { \delta_\sigma}(x_1^{\varepsilon}(i)){ \cdot}{{\bar{\sigma}}}(x^{\varepsilon}(i)),&\cdots(a)

\\[6pt]
    \mathbf{x}_2(\sigma(\varepsilon i))&=x_2^{ \mathrm{sign}({ {\sigma}}(\varepsilon i))}
(\varepsilon{ \delta_\sigma}(i){{\bar{\sigma}}}(i))\\[4pt]
    &=\varepsilon{ \delta_\sigma}(i){ \cdot}x_2^{\varepsilon{ \delta_\sigma}(i)}({{\bar{\sigma}}}(i)),&\cdots(b)

\\[6pt]
    \sigma(\mathbf{y}_1)(\varepsilon i)&=({ \delta_\sigma}{{\bar{\sigma}}})(\xi(y_1^{  \varepsilon_1},\varepsilon_1)(\varepsilon i)\cdot y_1^{\varepsilon_1}(\varepsilon i))&\\[4pt]
    &=\varepsilon \xi(y_1^{\varepsilon\varepsilon_1},\varepsilon_1)(i){ \delta_\sigma}(y_1^{\varepsilon\varepsilon_1(i)}(i))\cdot{{\bar{\sigma}}}(y_1^{\varepsilon\varepsilon_1(i)}(i)),&\cdots(c)

\\[6pt]
    \mathbf{y}_2(\sigma(\varepsilon i))&=\xi(y_2^{\varepsilon_2},\varepsilon_2)(\varepsilon{ \delta_\sigma}(i){{\bar{\sigma}}}(i))\cdot y_2^{\varepsilon_2}(\varepsilon{ \delta_\sigma}(i){{\bar{\sigma}}}(i)) &\\[4pt]
    &=\varepsilon{ \delta_\sigma}(i)\xi(y_2^{\varepsilon{ \delta_\sigma}\circ{{\bar{\sigma}}}^{-1}\varepsilon_2},\varepsilon_2)({{\bar{\sigma}}}(i))\cdot y_2^{\varepsilon{ \delta_\sigma}(i)\varepsilon_2({{\bar{\sigma}}}(i))}({{\bar{\sigma}}}(i)). &\cdots(d)
 \end{array}
\]
{ The right-hand side of each of $(a),(b),(c),(d)$ divides into the part in $\{\pm1 \}$ and the part in $I_d$ by dots, say the `$\{\pm1 \}$-part' and  the `$I_d$-part', respectively. 
The $\{\pm1 \}$-part of the equation $(a)=(b) $ implies (1). The $I_d$-part of the equation $(a)=(b) $ implies (2) by setting $i={\bar{\sigma}}^{-1}(i')$, $i'\in I_d$. }
For $\mathbf{y}_1$ and $\mathbf{y}_2$, { we have 
\begin{equation}
{{\bar{\sigma}}}(y_1^{\varepsilon\varepsilon_1(i)}(i))=y_2^{\varepsilon{ \delta_\sigma}(i)\varepsilon_2({{\bar{\sigma}}}(i))}({{\bar{\sigma}}}(i))
\end{equation}
 from the $I_d$-part of the equation $(c)=(d)$}. 
Setting {$\varepsilon :={ \delta_\sigma}(i)\varepsilon_2({{\bar{\sigma}}}(i))$ and $i={\bar{\sigma}}^{-1}(i')$, $i'\in I_d$}, we obtain (4). 
Similarly for 
{
the $\{\pm1 \}$-part of the equation $(d)=(c)$, we have
\[
 \left.
 \begin{array}{rlr}
\varepsilon{ \delta_\sigma}(i)\cdot \xi(y_2^{\varepsilon{ \delta_\sigma}\circ{{\bar{\sigma}}}^{-1}\varepsilon_2},\varepsilon_2)({{\bar{\sigma}}}(i))&=
\varepsilon \xi(y_1^{\varepsilon\varepsilon_1},\varepsilon_1)(i)\cdot{ \delta_\sigma}(y_1^{\varepsilon\varepsilon_1(i)}(i)),
\\[6pt]

{ \delta_\sigma}\circ{{\bar{\sigma}}}^{-1}(i')\cdot\xi(y_2,\varepsilon_2)(i')&= \xi(y_1^{\varepsilon\varepsilon_1},\varepsilon_1)({{\bar{\sigma}}}^{-1}(i'))\cdot{ \delta_\sigma}( y_1^{\varepsilon\varepsilon_1}({{\bar{\sigma}}}^{-1}(i')))&\\[4pt]

&=\xi({{\bar{\sigma}}}^\#y_1^{\varepsilon\varepsilon_1},\varepsilon_1\circ{{\bar{\sigma}}}^{-1})(i')\cdot{ \delta_\sigma}\circ{{\bar{\sigma}}}^{-1}({{\bar{\sigma}}}^\# y_1^{\varepsilon\varepsilon_1}(i'))&
\\[4pt]
&=\xi(y_2,\varepsilon_1\circ{{\bar{\sigma}}}^{-1})(i')\cdot{ \delta_\sigma}\circ{{\bar{\sigma}}}^{-1}(y_2(i')).  & \end{array}
 \right.
\]
So it follows that 
$\xi(y_2,\varepsilon_2)\cdot\xi(y_2,\varepsilon_1\circ{{\bar{\sigma}}}^{-1})\cdot \xi(y_2,{ \delta_\sigma}\circ{{\bar{\sigma}}}^{-1})
=1$. 
Using the formula $\xi(\tau,\lambda)\xi(\tau,\lambda')=\xi(\tau,\lambda\lambda')$}, we conclude (3). 

Suppose (1)-(4) conversely. Then for each $i\in{I}_d$, we have $(a)=(b)$ and $(c)=(d)$ for one of $\varepsilon\in\{\pm1 \}$. 
We may fill the equations for the other $\varepsilon\in\{\pm1 \}$ as follows.
First, the {$\{\pm1\}$-parts} of (a) and (b) coincide by (1).
The equality of the {$I_d$-parts} of (a) and (b) follows from (2) taking inverse mappings of both sides.
We can say the same for the {$I_d$-parts} of (c) and (d).
Finally, the equality of the {$\{\pm1\}$-parts} of (c) and (d) follows from (3) for each $y_2^{-1}(i)={{\bar{\sigma}}}^\#(y_1^{-\varepsilon_1\cdot\varepsilon_2\circ{{\bar{\sigma}}}\cdot{ \delta_\sigma}})(i)$, $i\in I_d$.
The induced translation map between $\mathcal{O}_1$ and $\mathcal{O}_2$ projects via their canonical double coverings by (1). 
\end{proof}

\begin{lemma}\label{calc_stratum}
Let $\mathcal{O}$ be an origami { of degree $d$} and $\pi_\mathcal{O}:\hat{\mathcal{O}}\rightarrow \mathcal{O}$ be the canonical double covering. 
Then, a point $p\in \mathcal{O}$ is a singularity of order $m$ $(m\in\{-1,0\}\cup\mathbb{N})$ if and only if either
\begin{enumerate}
\item  $m$ is even and $\pi_\mathcal{O}^{-1}(p)$ consists of two points of order $m$, or
 \item $m$ is odd and $\pi_\mathcal{O}^{-1}(p)$ consists of one point of order { $2m+2$}. 
\end{enumerate}
In particular, { 
if the permutation $\mathbf{z}_\mathcal{O}=\mathbf{x}_\mathcal{O}\mathbf{y}_\mathcal{O}\mathbf{x}_\mathcal{O}^{-1}\mathbf{y}_\mathcal{O}^{-1}$ has a cycle decomposition $a_1b_1\cdots a_kb_k c_1\cdots c_{k'}$ ($a_1b_1\cdots a_kb_k$, respectively) such that 
\begin{enumerate}
\item $a_j,b_j$ are cycles of length $m_j$ with $(i\mapsto-i)^\#(\mathbf{x}_\mathcal{O}\mathbf{y}_\mathcal{O})^\#a_j=b_j$, and  
\item $c_j$ is a cycle of length $m_j'$ with $(i\mapsto-i)^\#(\mathbf{x}_\mathcal{O}\mathbf{y}_\mathcal{O})^\#c_j=c_j$,  
\end{enumerate}
then the origami $\mathcal{O}$ belongs to the stratum $\mathcal{Q}_{\frac{d-k-k'}{2}}(2m_1-2,\ldots ,2m_k-2,m_1'-2,\ldots ,m_{k'}'-2)$ ($\mathcal{A}_{\frac{d-k}{2}}(2m_1-2,\ldots ,2m_k-2)$, respectively).
}
\end{lemma}
\begin{proof}
{
For the abelian origami $\hat{\mathcal{O}}=(\mathbf{x}_\mathcal{O},\mathbf{y}_\mathcal{O})$, each singularity of order $2m$ of the quadratic differential uniquely corresponds to a cycle of $\mathbf{z}_\mathcal{O}$ of length $m+1$. 
}The same holds for the half-rotated monodromy $(\mathbf{x}_\mathcal{O}\mathbf{y}_\mathcal{O})^\#\mathbf{z}_\mathcal{O}=\mathbf{x}_\mathcal{O}^{-1}\mathbf{y}_\mathcal{O}^{-1}\mathbf{x}_\mathcal{O}\mathbf{y}_\mathcal{O}$. 
{ By Remark \ref{canonical_double}, 
each two copies $p_1,p_2$ of a singularity $p$ of $\mathcal{O}$ of order $m$ are merged if $m$ is odd and separately left otherwise. 
The singularity $p$ has the cone angle $(m+2)\pi$ by (\ref{cone}). 
If  $p_1,p_2$ are merged, there are $4(m+2)$ squares summing up to a cycle of length $m+2$ around the resulting singularity on $\hat{\mathcal{O}}$. 
So the former claim follows. }

The involutive Deck transformation $\tau$ of $\pi_\mathcal{O}$ acts on $\bar{I}_d$ by $i\mapsto -i$ and on the monodromy group of $\hat{\mathcal{O}}$ by $\mathbf{x}_\mathcal{O}\mapsto \mathbf{x}_\mathcal{O}^{-1},\mathbf{y}_\mathcal{O}\mapsto \mathbf{y}_\mathcal{O}^{-1}$. 
{ Since $\tau\cdot\mathbf{z}_\mathcal{O}= (\mathbf{x}_\mathcal{O}\mathbf{y}_\mathcal{O})^\#\mathbf{z}_\mathcal{O}$, the correspondence of singularities of $\hat{\mathcal{O}}$ under $\tau$ can be seen as the $(\mathbf{x}_\mathcal{O}\mathbf{y}_\mathcal{O})$-conjugacy of the cycles of $\mathbf{z}_\mathcal{O}$}. 
With the former claim, we know that a singularity $p$ on $\mathcal{O}$ uniquely corresponds to either
\begin{enumerate}
\item two cycles $a,b$ of $\mathbf{z}_\mathcal{O}$ of length $m+1$ with $\tau\cdot a=(i\mapsto-i)^\#b$
 if $p$ has even order $2m$, or
\item one cycle $c$ of $\mathbf{z}_\mathcal{O}$ of length { $\frac{2m'+2}{2}+1=m'+2$} with $\tau\cdot c=(i\mapsto-i)^\#c$
 if $p$ has odd order $m'$. 
\end{enumerate}
So the latter claim follows. 
{ The genus of the origami $\mathcal{O}$ is given by the Euler charasteristic calculation. }
\end{proof}

\section{Algorithm}
\label{sec:4}
In this section, we present a series of algorithms for enumerating the equivalence classes of origamis of given degree $d$ and specifying their $PSL(2,\mathbb{Z})$-orbits. 
{ According to }
 Remark \ref{importantrem}, the { functioning} of these algorithms will be summarized as follows. 
\begin{theorem}\label{summary}
For each $d\in \mathbb{N}$, the combination of Algorithm { 4.2-4.7} outputs { generating systems and sets of representatives for} the Veech groups of all the equivalence classes of  origamis of degree $d$. 
\end{theorem}

A \textit{partition} \cite{AE} 
of $d$ is a finite sequence of weakly decreasing positive integers that sum to $d$. 
The \textit{partition number} $p(d)$, which counts the number of partitions of $d$, defines the following rapidly increasing sequence \cite{oeis}. 
\[
 \left.
 \begin{array}{l}
 1, 2, 3, 5, 7, 11, 15, 22, 30, 42, 56, 77, 101, 135, 176, 231, 297, 385, 490, 627, 792,\ldots 
 \end{array}
 \right.
\]
The following asymptotic formula \cite{HR} is known:
 \begin{equation}
    p(d)\sim\frac{1}{4d\sqrt{3}}\cdot e^{\sqrt{2d/3}}
 \end{equation}
The algorithm by Hashiguchi, Niki, and Nakagawa \cite{HNN} enumerates all the partitions of given integer.
We will accept the set $P(d)=\{(j_1,j_2,\ldots, j_n):\text{partition of }d\}$ as a known data.

To describe the class of each origami $\mathcal{O}=(x,y,\varepsilon)\in\bar{\Omega}_d$, we enumerate all the conjugators 
{
${\sigma}=\delta_\sigma\bar{\sigma}\in{ \bar{\mathfrak{S}}^{\mathrm{odd}}_{d}}$
}
satisfying the conditions in Lemma \ref{isom}. 
By (2) in Lemma \ref{isom}, up to equivalence, we only have to think of $x$ of the form
\begin{equation}\label{normalized_form}
  x=(1\cdots j_1)(j_1+1\cdots j_1+j_2)\cdots (\sum_{k=1}^{n-1} j_k+1\cdots d)
\end{equation} 
for each partition $(j_1,j_2,\ldots, j_n)\in P(d)$.
We will consider the \textit{restricted class} of an origami $\mathcal{O}$, the set of origamis with the same permutation `$x$' and equivalent to $\mathcal{O}$. 
By Lemma \ref{isom}, the restricted class is the conjugacy class in 
{
$\mathrm{Stab}(x):=\{{\sigma}=\delta_\sigma\bar{\sigma}\in{ \bar{\mathfrak{S}}^{\mathrm{odd}}_{d}}\mid\delta_\sigma =\delta_\sigma \circ x$ and $x={\bar{\sigma}}^\# (x^{\delta_\sigma }) \text{ on }I_d\}$.}
Remark that for general $y\in{ {\mathfrak{S}}_{d}}$ and $\varepsilon\in\mathcal{E}_{d}$, the mapping $y^{\varepsilon }:\bar{I}_d\rightarrow \bar{I}_d$ is not a bijection. 
It will be checked at (2) in Algorithm \ref{alg_isom}.

First, we present algorithms for enumerating all the $\mathrm{Stab}(x)$-conjugacy classes of origamis of the form $(x,y,\varepsilon)\in\bar{\Omega}_d$ satisfying the conditions in Lemma \ref{isom} for each $x\in P(d)$.
\begin{alg}\label{alg_isom}
  For each ${\mathcal{O}}=(x,y,\varepsilon)\in\bar{\Omega}_d$, we construct its restricted class $[{\mathcal{O}}]=\{(x,y',\varepsilon')\in\bar{\Omega}_d\mid (x,y',\varepsilon')\sim(x,y,\varepsilon)\}$ { where $\sim$ is the equivalence relation given by Lemma \ref{isom}} in the following steps: 
  \begin{enumerate}
    \item Take an element in $\mathrm{Stab}(x)$: ${\sigma}={\delta_\sigma}\bar{\sigma}\in{ \bar{\mathfrak{S}}^{\mathrm{odd}}_{d}}$ such that ${\delta_\sigma} ={\delta_\sigma} \circ x\text{ and }x={\bar{\sigma}}^\# (x^{{\delta_\sigma} }) \text{ on }I_d$.
    \item For each $\varepsilon'\in \mathcal{E}_d$, let $y_{{\sigma},\varepsilon'}:={\bar{\sigma}}^\#(y^{\varepsilon\cdot\varepsilon'\circ{\bar{\sigma}}\cdot{\delta_\sigma} })$. 
    Verify $\varepsilon'\in \mathcal{E}_d$ such that $y_{{\sigma},\varepsilon'}\in {\mathfrak{S}}_{d}$ and $\xi(y_{{\sigma},\varepsilon'},{\delta_\sigma} \circ{\bar{\sigma}}^{-1}\cdot \varepsilon\circ{\bar{\sigma}}^{-1}\cdot\varepsilon')=1$ on $I_d$.
    \item 
    Let $C_{{\sigma}}:=\{(x, y_{{\sigma},\varepsilon'},\varepsilon')\mid \varepsilon' $
    \ passes the test in (2)$\}$.
    \item Go back to (1) for some other leftover $\sigma\in\mathrm{Stab}(x)$. 
    When we have been through all elements in $\mathrm{Stab}(x)$, finish the algorithm and we conclude that $[\mathcal{O}]=\bigcup_{{\sigma}\in\mathrm{Stab}(x)}C_{{\sigma}}$.
  \end{enumerate}
\end{alg}
\begin{alg}\label{alg_class}
  Let $P(d)=\{(j_1,j_2,\ldots ,j_d):\text{partition of }d\}$. 
  We obtain the set $C\bar{\Omega}_d$ of the restricted classes of  all origamis of degree $d$ in the following steps.
    \begin{enumerate}
      \item $C\bar{\Omega}_d:=\emptyset$
      \item Take $j=(j_1,j_2,\ldots ,j_d)\in P(d)$. Define as follows:
      \[
       \left.
       \begin{array}{rl}
          d'_j&:=\mathrm{max}\{k\mid j_k>0\}, \\
          x_j&:=(1\ \cdots \ j_1)(j_1+1 \ \cdots \ j_1+j_2)\cdots(\sum_{k=1}^{d'-1} j_k+1\ \cdots \ d)\in\mathfrak{S}_d,
          \\R_j&:=\mathfrak{S}_d\times \mathcal{E}_d.
       \end{array}
       \right.
      \]
      \item Take $(y,\varepsilon)\in R_j$. 
      Apply Algorithm \ref{alg_isom} to $(x_j,y, \varepsilon)\in\bar{\Omega}_d$ to get $[(x_j,y, \varepsilon)]$.
      \item
      Add $[(x_j,y, \varepsilon)]$ to $C\bar{\Omega}_d$. 
      After that, remove $(y(\mathcal{O}),\varepsilon(\mathcal{O}))$ from $R_j$ for every $\mathcal{O}=(x_j,y(\mathcal{O}),\varepsilon(\mathcal{O}))\in [(x_j,y, \varepsilon)]$.
      \item Go back to (3) until $R_j=\emptyset$. If so, go to the next step.
      \item Go back to (2) for other leftover $j\in P(d)$. When we have been through all elements in $P(d)$, finish the algorithm.
    \end{enumerate}
\end{alg}


Next, 
{
we observe the action of $PSL(2,\mathbb{Z})$ on $\bar{\Omega}_d$ as shown in Remark \ref{importantrem} (2).  We calculate the two permutations $\varphi_T,\varphi_S\in\mathrm{Sym}(C\bar{\Omega}_d)$ defined by taking a new origami obtained by the decomposition into the pairs of directions $T(0,\frac{\pi}{2})=(0,\frac{\pi}{4})$ and $S(0,\frac{\pi}{2})=(-\frac{\pi}{2},0)$, respectively. 
}

Let $C\bar{\Omega}_d$ be the output in Algorithm \ref{alg_class}. 

\begin{figure}[htbp]
\begin{center}
  \includegraphics[width=130mm]{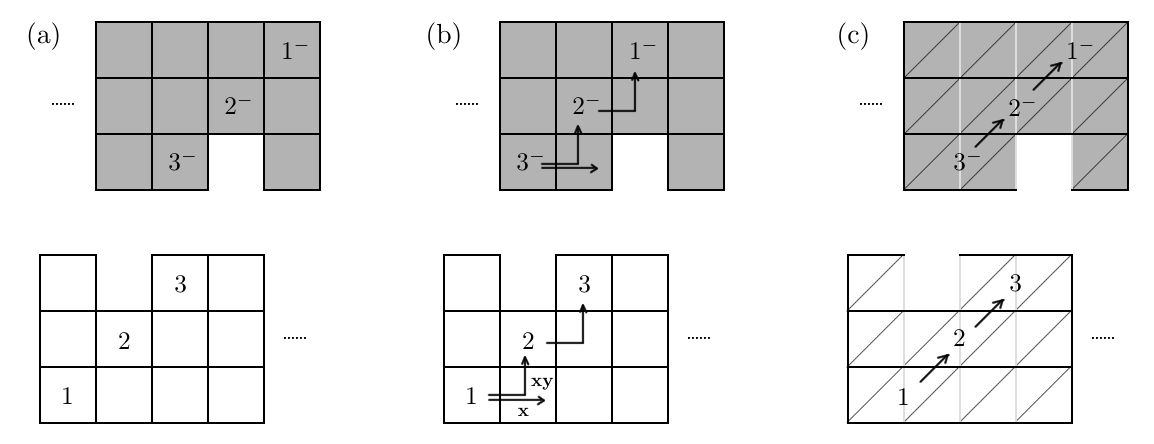}
  \caption{Decomposition of origami (a) into $(0,\frac{\pi}{4})$ { with { a} double symmetry}: 
The desired decomposition (c) is obtained { as (b): applying $\gamma_T$ and taking the conjugation in $(-i \mapsto -x^{-1}(i)\mid i \in I_d)$. }
}
\label{Figure3}
\end{center}
\end{figure}

{\textbf{1.\ }}
To obtain the permutation $\varphi_T$, we consider as follows:
\[
(x,y,\varepsilon)\xmapsto{\text{Def.\ref{xye_def}}}(\mathbf{x},\mathbf{y})\overset{\gamma_T}{\longmapsto}\overset{\text{conj.}}{\longmapsto}(\mathbf{x}_T,\mathbf{y}_T)\xmapsto{\text{\text{Lem.\ref{restore}}}}(x_T,y_T,\varepsilon_T). 
\]
To apply Lemma \ref{restore}, 
{
we have to find a conjugator yealding { a} double symmetry of $\mathbf{x}_T$, a cycle decomposition of $\mathbf{y}_T$, and $\varepsilon_{\mathbf{y}_T}$.
}
Remark that the decomposition into $T(0,\frac{\pi}{2})=(0,\frac{\pi}{4})$ { with { a} double symmetry} is given by $\gamma_T$ and the conjugation in $(-i \mapsto { -}x^{-1}(i)\mid i \in I_d)$ as shown in Figure\ \ref{Figure3}.

For $\mathcal{O}=(x,y,\varepsilon)\in \bar{\Omega}_d$, $a\in I_d$, and $\varepsilon'\in\{\pm 1\}$,  we have:
\begin{align}
  \gamma_T(\mathbf{y}_\mathcal{O})(\varepsilon' a)
  &=\mathbf{y}_\mathcal{O}\circ\mathbf{x}_\mathcal{O}(\varepsilon' a)
  \nonumber\\
  &=\mathbf{y}_\mathcal{O}(\varepsilon' x^{\varepsilon'}(a))
  \nonumber\\
  &=\varepsilon(\varepsilon' x^{\varepsilon'}(a)))\cdot y^{\varepsilon(\varepsilon' x^{\varepsilon'}(a)))}(\varepsilon' x^{\varepsilon'}(a)))\cdot\varepsilon(y^{\varepsilon(\varepsilon' x^{\varepsilon'}(a)))}(\varepsilon' x^{\varepsilon'}(a)))
  \nonumber\\
  &=\varepsilon'\varepsilon( x^{\varepsilon'}(a)))\cdot \varepsilon'y^{\varepsilon'\varepsilon( x^{\varepsilon'}(a)))}( x^{\varepsilon'}(a))\cdot \varepsilon'\varepsilon(y^{\varepsilon'\varepsilon( x^{\varepsilon'}(a)))}( x^{\varepsilon'}(a)))
  \nonumber\\
  &=\varepsilon'\varepsilon( x^{\varepsilon'}(a)))\cdot\varepsilon(y^{\varepsilon'\varepsilon( x^{\varepsilon'}(a)))}( x^{\varepsilon'}(a)))\cdot y^{\varepsilon'\varepsilon( x^{\varepsilon'}(a)))}( x^{\varepsilon'}(a)).
  \label{gammaT_calc}
\end{align}

\begin{alg}\label{alg_T}
  Let $C=[(x,y,\varepsilon)]\in C\bar{\Omega}_d$ be a restricted class. 
  By (\ref{gammaT_calc}), we obtain $\varphi_T(C)$ in the following steps:
  \begin{enumerate}
    \item $I_d':=I_d$, $j:=0$.
    \item $a_{0,j}:=\mathrm{min}(I_d')$, $\varepsilon_{0,j}':=1$, $i:=0$.
    \item $b_{i,j}:=x^{\varepsilon_{i,j}'}(a_{i,j})$,  $a_{i+1,j}':=y^{\varepsilon_{i,j}'\varepsilon(b_{i,j})}(b_{i,j})$,  $\varepsilon_{i+1,j}':=\varepsilon_{i,j}'\varepsilon(b_{i,j})\varepsilon({ a_{i+1,j}'})$.
    \item Remove $a_{i,j}$ from $I_d'$.
    Define:
    \[
    a_{i+1,j}:= \left\{
     \begin{array}{rl}
        a_{i+1,j}'&\text{ if $\varepsilon_{i+1,j}'=1$,}\\
        x^{-1}(a_{i+1.j}')&\text{otherwise.}
     \end{array}
     \right. 
    \]
    
    \item If $a_{i+1,j}=a_{0,j}$, let $c_j:=(a_{0,j}\ a_{1,j}\ \cdots \  a_{i,j})$. 
    Otherwise, go back to (3) for the next $i$.
    \item If $I_d'\neq \emptyset $, go back to (2) for the next $j$. 
    Otherwise, finish the loop and let  $x_T=x$, $y_T:=c_1c_2\cdots c_j$, and $\varepsilon_T:=(a_{i,j}\mapsto\varepsilon_{i,j}')$.
    \item { Search} for the equivalence class $C_T\in C\bar{\Omega}_d$ represented by $(x_T,y_T,\varepsilon_T)$ and we conclude that $\varphi_T(C)=C_T$.
  \end{enumerate}
\end{alg}

{\textbf{2.\ }}
To obtain the permutation $\varphi_S$, we consider as follows:
\[
(x,y,\varepsilon)\xmapsto{\text{Def.\ref{xye_def}}}(\mathbf{x},\mathbf{y})\overset{\gamma_S}{\longmapsto}
\overset{\text{conj.}}{\longmapsto}(\mathbf{x}_S,\mathbf{y}_S)\xmapsto{\text{Lem.\ref{restore}}}\overset{\text{conj.}}{\longmapsto}(x_S,y_S,\varepsilon_S).
\]
We use two conjugators in $\bar{\mathfrak{S}}_d$: the former collects signs of cells in each vertical cylinder to apply Lemma \ref{restore}, and the latter makes $x_S$ to be the normalized form (\ref{normalized_form}).
The former conjugator is given by $\sigma_\delta:=(\pm i\mapsto \pm\delta(i)i\mid i\in I_d)\in{ \bar{\mathfrak{S}}^\mathrm{odd}_d}$ where $\delta\in\mathcal{E}_d$ satisfies that for every cycle $c$ in $\mathbf{x}$, $\{\delta(|i|)i\mid i\in c\}$ forms a cycle either $c$ or $c'$.

For $\mathcal{O}=(x,y,\varepsilon)\in \bar{\Omega}_d$, $a\in I_d$, and $\delta'\in\{\pm 1\}$, we have:
\begin{align}
  \gamma_S(\mathbf{x})(\delta' a)
  &=\mathbf{y}(\delta' a)
  \nonumber\\
  &=\varepsilon(\delta' a)\cdot y^{\varepsilon(\delta' a)}(\delta' a)\cdot\varepsilon(y^{\varepsilon(\delta' a)}(\delta' a))
  \nonumber\\
  &=\delta'\varepsilon( a)\cdot \delta' y^{\delta'\varepsilon( a)}( a)\cdot\delta' \varepsilon(y^{\delta' \varepsilon(a)}(a))
  \nonumber\\
  &=\delta'\varepsilon( a)\varepsilon(y^{\delta' \varepsilon(a)}(a))\cdot y^{\delta'\varepsilon( a)}( a).
  \label{delta_calc}
\end{align}
\begin{alg}\label{alg_delta}
  Let $C=[(x,y,\varepsilon)]\in C\bar{\Omega}_d$ be a restricted class. By (\ref{delta_calc}),  we obtain the required $\delta=\delta(C)\in\mathcal{E}_d$ in the following steps:
  \begin{enumerate}
  \item $I_d':=I_d$, $j:=0$.
  \item $a_{0,j}:=\mathrm{min}(I_d')$, $\delta_{0,j}:=1$, $i:=0$.
  \item $a_{i+1,j}:=y^{\delta_{i,j} \varepsilon(a_{i,j})}(a_{i,j})$, $\delta_{i+1,j}:=\delta_{i,j}\varepsilon( a_{i,j})\varepsilon(a_{i+1,j})$
  \item Remove $a_{i,j}$ from $I_d'$.
  \item If $a_{i+1,j}=a_{0,j}$, let $c_j:=(a_{0,j}\  a_{1,j}\ \ldots \  a_{i,j})$. Otherwise go back to (3) for the next $i$.
  \item If $I_d'\neq \emptyset $ then go back to (2) for the next $j$. Otherwise finish the loop and let  $x_S':=c_1c_2\cdots c_j$ and $\delta:=(a_{i,j}\mapsto\delta_{i,j})$.
  \end{enumerate}
\end{alg}

To apply Lemma \ref{restore}, we will calculate $\varepsilon_{{\sigma_\delta}^\#\mathbf{y}_S}$ and a cycle decomposition of ${\sigma_\delta}^\#\mathbf{y}_S$.
After that, we apply the conjugator which makes $x_S$ to the normalized form (\ref{normalized_form}). 
So in advance, we will prepare the list $\{\sigma^\#x_p\mid \sigma\in\mathfrak{S}_d\}$ equipped with information of conjugator for each $p\in P(d)$. 
Note that the permutations `$x$' of any two equivalent origamis share the same partition by Lemma \ref{isom}. 
So the restricted classes calculated from Algorithm \ref{alg_isom} with this list exhausts all the patterns of origamis.

For $(x,y,\varepsilon)\in \bar{\Omega}_d$, $a\in I_d$ and $\varepsilon'\in\{\pm 1\}$, we have:
\begin{align}
  { {\sigma_\delta}^\#}(\mathbf{y}_S)(\varepsilon' a)
  &={ {\sigma_\delta}}(\mathbf{x}^{-1}(\delta(|\varepsilon' a|)\varepsilon' a))
  \nonumber\\
  &=\delta(|\mathbf{x}^{-1}(\varepsilon'\delta(a) a)|)\cdot \mathbf{x}^{-1}(\varepsilon'\delta(a) a)
  \nonumber\\
  &=\varepsilon'\delta(a)\cdot \delta(x^{-\varepsilon'\delta(a)}(a))\cdot x^{-\varepsilon'\delta(a)}(a).
  \label{gammaS_calc}
\end{align}

\begin{alg}\label{alg_S}
  Let $C=[(x,y,\varepsilon)]\in C\bar{\Omega}_d$ be a restricted class and $\delta=\delta(C)\in\mathcal{E}_d$ be the output in Algorithm \ref{alg_delta}. 
  By (\ref{gammaS_calc}), we obtain $\varphi_S(C)$ in the following steps:
  \begin{enumerate}
    \item $I_d':=I_d$, $j:=0$.
    \item $a_{0,j}:=\mathrm{min}(I_d')$, $\varepsilon_{0,j}':=1$, $i:=0$
    \item Remove $a_{i,j}$ from $I_d'$.
    Let $a_{i+1,j}:=x^{-\varepsilon_{i,j}'\delta(a_{i,j})}(a_{i,j})$, $\varepsilon_{i+1,j}':=\varepsilon_{i,j}'\delta(a_{i,j})\delta(a_{i+1,j})$. 
    \item If $a_{i+1,j}=a_{0,j}$, let $c_j:=(a_{0,j}\ a_{1,j}\ \cdots \ a_{i,j})$. Otherwise go back to (3) for the next $i$.
    \item If $I_d'\neq \emptyset $, go back to (2) for the next $j$. 
    Otherwise finish the loop and let $x_S':=\delta^\#x_S$, $y_S':=c_1c_2\cdots c_j$ and $\varepsilon_S':=(a_{i,j}\mapsto\varepsilon_{i,j}')$.
    \item { Search} for the conjugator $\sigma\in\mathfrak{S}_d$ such that $\sigma^\#x_S'$ is of normalized form. 
    Let $(x_S,y_S,\varepsilon_S):=(\sigma^\#x_S',\sigma^\#y_S',\varepsilon_S'\circ\sigma^{-1})$.
    \item { Search} for the equivalence class $C_S\in C\bar{\Omega}_d$ represented by $(x_S,y_S,\varepsilon_S)$ and we conclude that $\varphi_S(C)=C_S$.
  
  \end{enumerate}
\end{alg}


\begin{alg}\label{alg_comp}
Let $\varphi_T,\varphi_S\in\mathrm{Sym}(C\bar{\Omega}_d)$. 
We obtain the $\langle\varphi_T^{-1},\varphi_S^{-1}\rangle$-orbit decomposition $C\bar{\Omega}_d{ =\sqcup_{t=1,2,...} O_t}$ in the following steps.
  \begin{enumerate}\setlength{\itemsep}{4pt} 
    \item { $R:=C\bar{\Omega}_d$, $t:=1$.}
    \item { $O_t:=\emptyset$.}
    \item { Take $C\in R$ and add $C$ to $O_t$.}
    \item { Take $C\in O_t$ and let $O(C):=\{\varphi_T^{-k}(C),\varphi_S^{-k}(C)\mid k\in\mathbb{N} \}$. }
    \item Add all elements in ${ O(C)}$ to $O_t$ and remove them from ${ R}$.
    \item Go back to (4) for other leftover ${ C}\in O_t$. 
When we have been through all elements in $O_t$, go to the next step.
    \item Go back to (2) for the next $t$ until ${ R}=\emptyset$. If so, finish the algorithm.
  \end{enumerate}
\end{alg}
{ For each orbit $O_t\subset C\bar{\Omega}_d$ in the output and an origami $[\mathcal{O}]\in O_t$, the universal Veech group $PSL(2,\mathbb{Z})=\langle T,S\rangle$ acts on $O_t$ by $T\cdot [\mathcal{O}]=\varphi_T^{-1}([\mathcal{O}])$ and $S\cdot [\mathcal{O}]=\varphi_S^{-1}([\mathcal{O}])$ (see also Remark \ref{importantrem}). 
The Teichm\"uller curve $C(\mathcal{O})$ corresponds to the orbit $O_t$, and the Veech group  is described by 
\begin{equation}\Gamma(\mathcal{O})=\{w(T,S)\in PSL(2,\mathbb{Z})\mid w:\text{word,\ } w(\varphi_T^{-1},\varphi_S^{-1}) [\mathcal{O}]=[\mathcal{O}]\}.\end{equation}
A generating system and a set of representatives is calculated by the Reidemeister-Schreier method \cite{MKS}. } 

\section{Calculation results}
\label{sec:5}
In the following, we show some calculation results obtained by the algorithms stated in the previous section powered by Python.
We first note that all classes representing disconnected origamis are removed and abelian origamis are treated as flat surfaces. 

\begin{figure}[p]
\begin{center}
  \includegraphics[width=120mm]{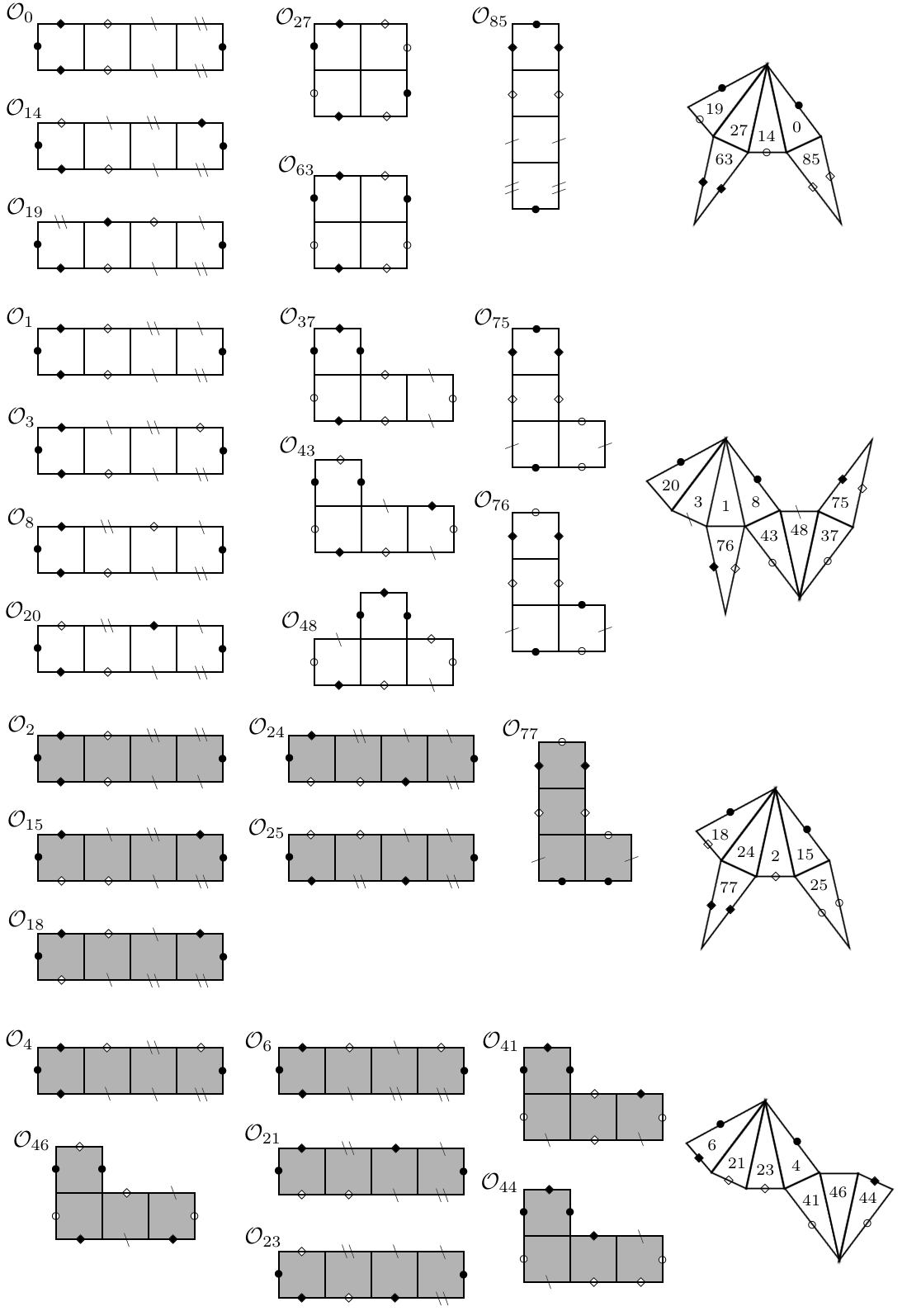}
\caption{(Part 1/2) All classes of origamis of degree $4$ and their positions in Teichm\"uller curves. { Non-abelian origamis are shaded. }}
\label{FigureA}
\end{center}
\end{figure}
\begin{figure}[p]
\begin{center}
  \includegraphics[width=120mm]{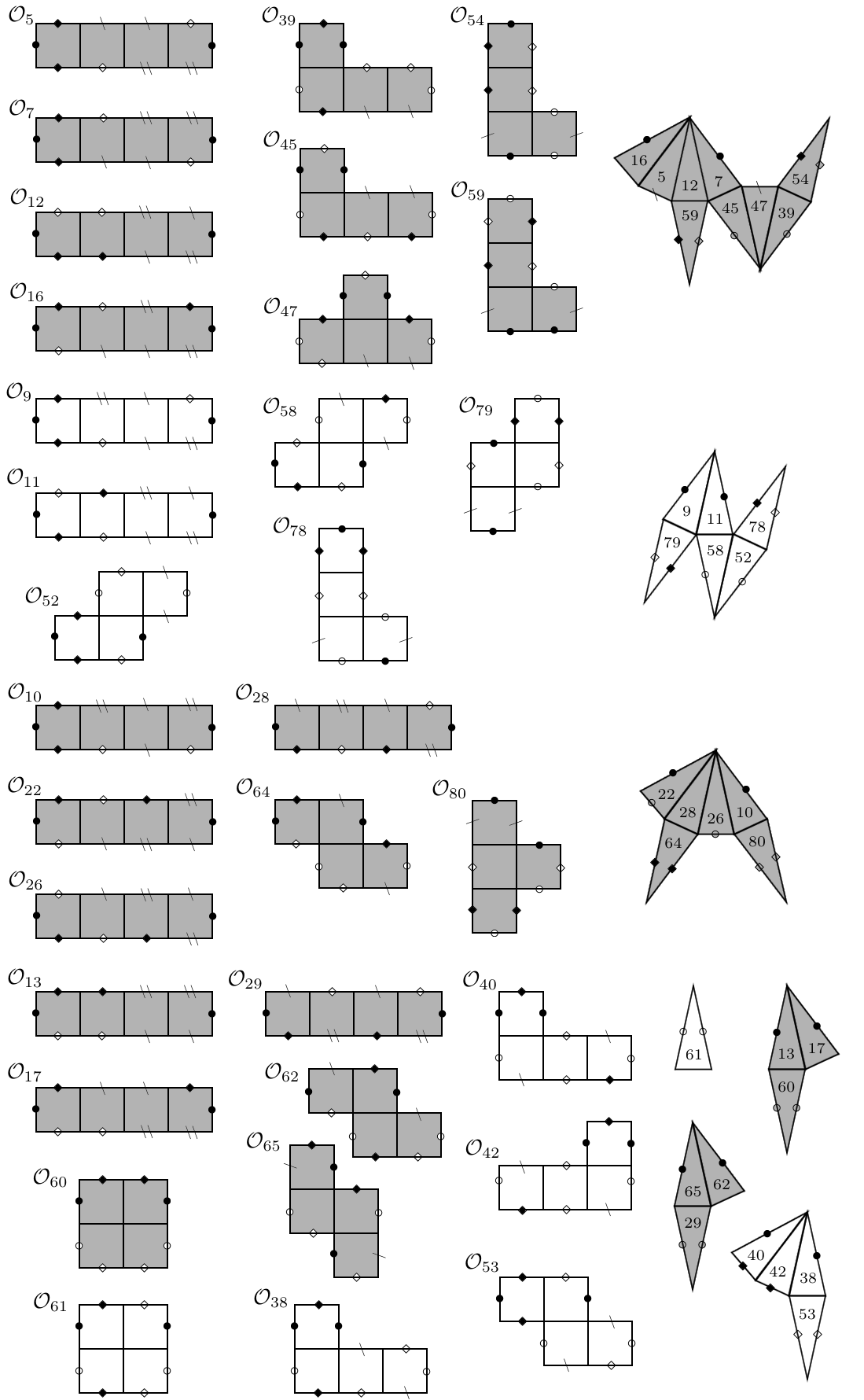}
\caption{(Part 2/2) All classes of origamis of degree $4$ and their positions in Teichm\"uller curves. { Non-abelian origamis are shaded. }}
\label{FigureB}
\end{center}
\end{figure}

Figure \ref{FigureA} and Figure \ref{FigureB} show all classes of origamis of degree $4$ and their positions in Teichm\"uller curves. The classes of origamis are numbered according to Algorithm \ref{alg_class} 
{ with the following orders: the graded reverse {lexicographic order} \cite{AE} on $P(d)$ graded by the number of parts, the lexographic order with respect to the vector representation $(y(1),...,y(d))$, $y\in \mathfrak{S}_d$, and the byte order of $\mathcal{E}_d$. 
Concrete data on this numbering of origamis can be found in \cite{GitK}.}
There are $26$ classes of abelian origamis summing up to $5$ components and $34$ classes of non-abelian origamis summing up to $6$ components.
In the figures, we denote the copies of the standard fundamental domain of the group $PSL(2,\mathbb{Z})$ by isosceles triangles where the keen vertices correspond to the cusps.
Every two edges with the same symbol are glued so that the cusps match.
Every edge with no symbol is glued individually, making a conical point of angle $\pi/2$.


Table \ref{Table1} shows the number of classes of origamis, the number of components of Teichm\"uller curves, the range of genus of Teichm\"uller curves, and the number of classes of possible $G_\mathbb{Q}$-conjugacy for degree $1\leq d\leq7$. 
The possibility of $G_\mathbb{Q}$-conjugacy is checked by the following $G_\mathbb{Q}$-invariants: 
the degree and the stratum of an origami $\mathcal{O}$; the degree 
, the genus, and the valency list of $C(\mathcal{O})$. 
We calculate strata of origamis by Lemma \ref{calc_stratum}. 
The distribution of the number of classes origamis per genus is shown in Table \ref{genus}.
\begin{table}[H]

\begin{center}
\scalebox{0.65}{
\begin{tabular}{|c||c|c|c|c||c|c|c|c|} \hline
  & \multicolumn{4}{|c||}{abelian}& \multicolumn{4}{|c|}{non-abelian} \\ \hline
  $d$& \#$\{\mathcal{O}\}$ & \#$\{C(\mathcal{O})\}$ &$g(C(\mathcal{O}))$&possible $G_\mathbb{Q}$-conjugacy& \#$\{\mathcal{O}\}$ & \#$\{C(\mathcal{O})\}$&$g(C(\mathcal{O}))$  &possible $G_\mathbb{Q}$-conjugacy\\ \hline
  1& 1 &1&0&none&0&0&0&none\\
  2& 3 &1&0&none&1&1&0&none\\
  3& 7 &2&0&none&4&1&0&none\\
  4& 26 &5&0&none&34&6&0&none\\
  5& 91 &8&0&none&227&13&0&none\\
  6& 490 &28&0&1 class&2316&88&0&13 classes\\
  7& 2785 &41&$0\sim1$&5 classes&26574&88&$0\sim11$ &3 classes\\\hline
\end{tabular}
}
\caption{Summary of the result for degree $1\leq d\leq7$}\label{Table1}
\end{center}
\end{table}
\begin{table}[htb]
\begin{center}
\scalebox{0.65}{
\begin{tabular}{|c||c|c|c|c|c|c||c|c|c|c|c|c|} \hline
  & \multicolumn{6}{|c||}{abelian}& \multicolumn{6}{|c|}{non-abelian} \\ \hline
  $d\setminus g$& 0&1&2&3&4 &total& 0&1&2&3&4 &total\\ \hline
  $1$                  & 0&1&0&0&0 &1& 0&0&0&0&0 &0\\
  $2$                  & 0&3&0&0&0 &3& 1&0&0&0&0 &1\\
  $3$                  & 0&4&3&0&0 &4& 0&4&0&0&0 &4\\
  $4$                  & 0&7&19&0&0 &26& 3&15&16&0&0 &34\\
  $5$                  & 0&6&51&34&0 &91& 0&49&138&40&0 &227\\
  $6$                  & 0&15&142&333&0 &490& 10&181&1085&1040&0 &2316\\
  $7$                  & 0&8&250&1735&792 &2785& 0&534&6449&16000&3591 &26574\\\hline
\end{tabular}
}
\caption{The number of  classes of origamis of given genus and degree.}\label{genus}
\end{center}
\end{table}

\begin{figure}[htbp]
\begin{center}
  \includegraphics[width=125mm]{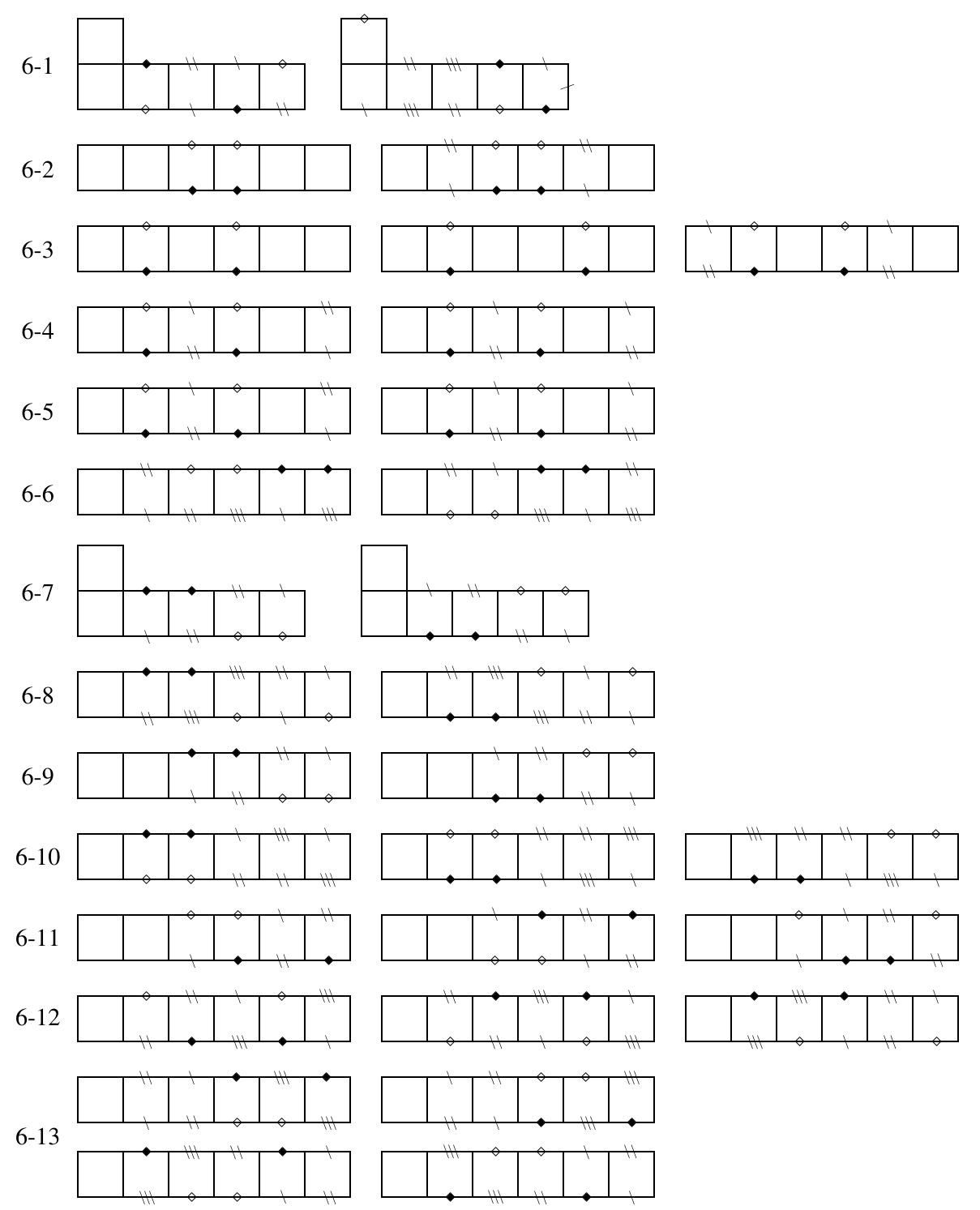}
  \caption{Origamis that induce Teichm\"uller curves in Table \ref{Table2}: unmarked edges are glued with the opposite. }
\label{FigureD}
\end{center}
\end{figure}

\begin{figure}[htbp]
\begin{center}
  \includegraphics[width=95mm]{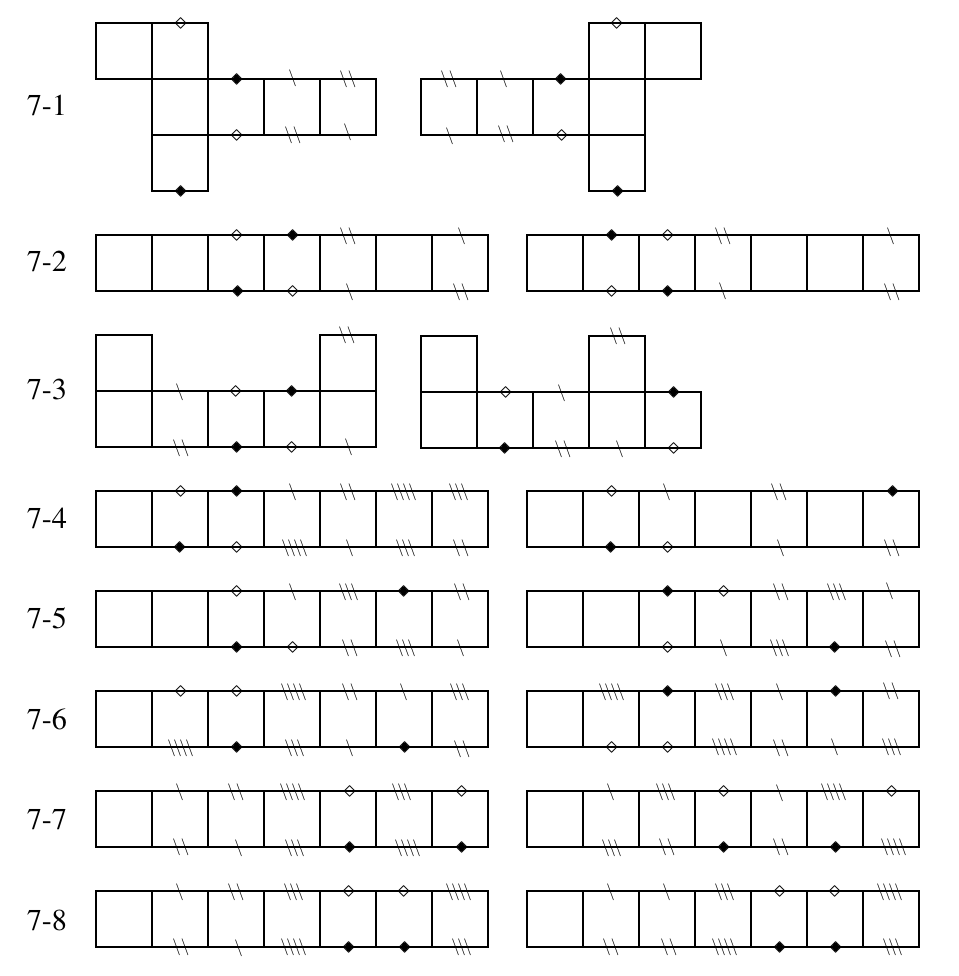}
  \caption{Origamis that induce Teichm\"uller curves in Table \ref{Table3}: unmarked edges are glued with the opposite. }
\label{FigureE}
\end{center}
\end{figure}

The more detailed description of classes of possible $G_\mathbb{Q}$-conjugacy for degree $d\leq7$ is as follows.
\begin{theorem}
All the Teichm\"uller curves induced from origamis of degree $d\leq7$ are distinguished by Galois invariants except for the 13 cases in Table \ref{Table2} and the 8 cases in Table \ref{Table3}. 
Figure \ref{FigureD} and Figure \ref{FigureE} show origamis that induce Teichm\"uller curves in each of the exceptional cases.  
\end{theorem}

\begin{table}[htb]
\begin{center}
\scalebox{0.65}[0.7]{
\begin{tabular}{|c|c|c|c|c|} \hline
  No. & stratum  &index& valency list of $C(\mathcal{O})$& relationship between $C(\mathcal{O})$   \\ \hline
  6-1 &$\mathcal{A}_3(0,8)$&15&$(3^5\mid2^7,1\mid5,4,3^2)$&two identical, mirror-closed curves\\
  6-2 &$\mathcal{Q}_1(-1^2,0^3,2)$&12&$(3^4\mid 2^6\mid6,3,2,1)$&two identical, mirror-closed curves\\
  6-3 &$\mathcal{Q}_2(-1^2,0,6)$&12&$(3^4\mid 2^6\mid6,3,2,1)$&two identical, mirror-closed curves\\
  6-4 &$\mathcal{Q}_2(0^2,2^2)$&12&$(3^4\mid 2^6\mid6,3,2,1)$&three identical, mirror-closed curves\\
  6-5 &$\mathcal{Q}_3(2,6)$&12&$(3^4\mid 2^6\mid6,3,2,1)$&two identical, mirror-closed curves\\
  6-6 &$\mathcal{Q}_2(-1^2,3^2)$&15&$(3^5\mid 2^7,1\mid6,5,3,1)$&two distinct, mirror-closed curves\\
  6-7 &$\mathcal{Q}_2(-1^2,3^2)$&15&$(3^5\mid 2^7,1\mid5,4,3^2)$&one pair of mirror-symmetric curves, mirroring each other\\
  6-8 &$\mathcal{Q}_3(-1,9)$&22&$(3^7,1 \mid 2^{11}\mid6,5,4^2,3)$&one pair of mirror-symmetric curves, mirroring each other\\
  6-9 &$\mathcal{Q}_2(-1^2,0,6)$&24&$(3^8\mid 2^{12}\mid6,5,4^2,3,2)$&one mirror-conjugate pair\\
  6-10 &$\mathcal{Q}_2(-1^3,7)$&27&$(3^9\mid 2^{13},1\mid6^2,5,4,3^2)$&one mirror-conjugate pair \& one mirror-closed curve\\
  6-11 &$\mathcal{Q}_2(-1,0,1,4)$&36&$(3^{12}\mid 2^{18}\mid6^2,5^2,4^2,3^2)$&one mirror-conjugate pair \& one mirror-closed curve\\
  6-12 &$\mathcal{Q}_3(1,7)$&54&$(3^{18}\mid 2^{27}\mid6^4,5^3,4^3,3)$&one mirror-conjugate pair \& one mirror-closed curve\\
  6-13 &$\mathcal{Q}_3(-1,9)$&66&$(3^{22}\mid 2^{33}\mid6^6,5^3,4^3,3)$&two mirror-conjugate pairs\\\hline
\end{tabular}
}
\caption{Classes of possible $G_\mathbb{Q}$-conjugacy for degree $6$}\label{Table2}
\end{center}
\end{table}

\begin{table}[htb]
\begin{center}
\scalebox{0.65}[0.7]{
\begin{tabular}{|c|c|c|c|c|} \hline
  No. & stratum  &index& valency list of $C(\mathcal{O})$& relationship between $C(\mathcal{O})$   \\ \hline
  7-1 &$\mathcal{A}_4(12)$&7&$(3^2,1\mid2^3,1\mid4, 3)$&one pair of mirror-symmetric curves, mirroring each other\\
  7-2 &$\mathcal{A}_3(0,2,6)$&16&$(3^5,1\mid 2^8\mid7,4,3,2)$&two distinct, mirror-closed curves\\
  7-3 &$\mathcal{A}_4(12)$&21&$(3^7\mid 2^{11}\mid6,5,4,3^2)$&two distinct, mirror-closed curves\\
  7-4 &$\mathcal{A}_4(12)$&42&$(3^{14}\mid 2^{21}\mid7^2,5^2,4^3,3^2)$&two distinct, mirror-closed curves\\
  7-5 &$\mathcal{A}_3(0,2,6)$&48&$(3^{16}\mid 2^{24}\mid7^2,6,5^2,4^3,3^2)$&one mirror-conjugate pair\\
  7-6 &$\mathcal{Q}_2(-1,1^3,2)$&16&$(3^5,1\mid 2^8\mid7,6,2,1)$&one mirror-conjugate pair\\
  7-7 &$\mathcal{Q}_4(12)$&28&$(3^9,1\mid 2^{14}\mid7^2,6,3^2,2)$&two distinct, mirror-closed curves\\
  7-8 &$\mathcal{Q}_3(-1^2,10)$&36&$(3^{12}\mid 2^{18}\mid7^3,6,3^2,2,1)$&two distinct, mirror-closed curves\\\hline
\end{tabular}
}
\caption{Classes of possible $G_\mathbb{Q}$-conjugacy for degree $7$}\label{Table3}
\end{center}
\end{table}

\begin{remark}
The mirror relation between surfaces implies a Galois conjugacy which induces complex conjugacy. 
For example, the situation `one pair of mirror-symmetric curves, mirroring each other' in Table \ref{Table2}, \ref{Table3} is caused by such a Galois conjugacy modifying the embeddings of Teichm\"uller curves into the moduli space. 
If there is no mirror relation, we only know the invariant agreement and do not know if there is a nontrivial $G_\mathbb{Q}$-conjugacy.
\end{remark}

\newpage
\begin{DAS}
The data that support the findings of this study are openly available in GitHub \cite{GitK}. 
\end{DAS}
\begin{acknowledgments}
{ The author} thanks to Prof.\ Toshiyuki Sugawa and Prof.\ Hiroshige Shiga for their helpful advices and comments. 
{ The author thanks the anonymous referee for his careful reading of the manuscript and his suggestions for improvement. }
The computation was carried out using the computer resource offered under the category of General Projects by Cyber Science Center, Tohoku University.
\end{acknowledgments}

%
%

\bibliographystyle{abbrv}     
\bibliography{references.bib}   

\begin{description}
\item \ \\Shun Kumagai\\Research Center for Pure and Applied Mathematics \\Graduate School of Information Sciences \\Tohoku University \\Sendai, 980-8579 \\Japan \\Email: shun.kumagai.p5@alumni.tohoku.ac.jp

\end{description} 

\end{document}